\documentclass[10pt]{article}
\usepackage{amssymb}
\usepackage{amsmath}
\usepackage{amscd}
\newtheorem{theorem}{Theorem}[section]
\newtheorem{proposition}[theorem]{Proposition}
\newtheorem{lemma}[theorem]{Lemma}
\newtheorem{corollary}[theorem]{Corollary}

\newtheorem{condition}[theorem]{Condition}
\newtheorem{definition}[theorem]{Definition}
\newtheorem{remark}[theorem]{Remark}
\newtheorem{example}[theorem]{Example}

\newcommand{\bthm}{\begin{theorem}}
\newcommand{\ethm}{\end{theorem}}
\newcommand{\bpr}{\begin{proposition}}
\newcommand{\epr}{\end{proposition}}
\newcommand{\blem}{\begin{lemma}}
\newcommand{\elem}{\end{lemma}}
\newcommand{\bco}{\begin{corollary}}
\newcommand{\eco}{\end{corollary}}
\newcommand{\bde}{\begin{definition}\rm}
\newcommand{\ede}{\end{definition}}
\newcommand{\bre}{\begin{remark}\rm}
\newcommand{\ere}{\end{remark}}
\newcommand{\bex}{\begin{example}\rm}
\newcommand{\eex}{\end{example}}
\newcommand{\bcon}{\begin{condition}\rm}
\newcommand{\econ}{\end{condition}}
\newcommand{\bprf}{\noindent{\it Proof.\ }}
\newcommand{\eprf}{\hspace*{\fill} \rule{1.6mm}{3.2mm} \vspace{1.8mm}}
\newcommand{\benu}{\begin{enumerate}\renewcommand{\labelenumi}{{\rm (\roman{enumi})}}\renewcommand{\itemsep}{0pt}}
\newcommand{\eenu}{\end{enumerate}}

\newcommand{\N}{\mathbb{N}}
\newcommand{\Z}{\mathbb{Z}}
\newcommand{\R}{\mathbb{R}}
\newcommand{\C}{\mathbb{C}}
\newcommand{\T}{\mathbb{T}}
\newcommand{\e}{\varepsilon}
\newcommand{\K}{\mathbb{K}}
\newcommand{\M}{\mathbb{M}}
\newcommand{\I}{\mathbb{I}}

\DeclareMathOperator{\Aut}{Aut}
\DeclareMathOperator{\Prim}{Prim}
\DeclareMathOperator{\spa}{span}
\DeclareMathOperator{\cspa}{\overline{span}}
\newcommand{\On}{{\mathcal O}_n}

\newcommand{\cp}{{\mathcal O}_n{\rtimes_{\alpha^\omega}}G}

\newcommand{\cpr}{{\mathcal O}_n{\rtimes_{\alpha^\omega}}\mathbb{R}}

\newcommand{\ip}[2]{\langle\,{#1}\,|\,{#2}\,\rangle}
\newcommand{\W}{{\mathcal W}}

\begin{document}
\title{The ideal structures of crossed products of Cuntz algebras by quasi-free actions of abelian groups}
\author{Takeshi KATSURA\\
Department of Mathematical Sciences\\
University of Tokyo, Komaba, Tokyo, 153-8914, JAPAN\\
e-mail: {\tt katsu@ms.u-tokyo.ac.jp}}
\date{}

\maketitle

\begin{abstract}
{\footnotesize We completely determine the ideal structures of the crossed products of Cuntz algebras by quasi-free actions of abelian groups and give another proof of A. Kishimoto's result on the simplicity of such crossed products. We also give a necessary and sufficient condition that our algebras become primitive, and compute the Connes spectra and K-groups of our algebras.}
\end{abstract}

\section{Introduction}\label{INTRO}

Recently the classification theory of simple $C^*$-algebras has developed rapidly.
One of the most important questions in the classification theory of $C^*$-algebras is to determine whether a given $C^*$-algebra is simple or not.
It is also important to examine the ideal structure of a given $C^*$-algebra if it turns out to be non-simple.
There have been many works examining the ideal structures of some classes of $C^*$-algebras.
J. Cuntz examined the ideal structures of Cuntz-Krieger algebras 
under a certain condition in \cite{C2}.
In \cite{aHR}, A. an Huef and I. Raeburn determined the ideal structures
 of arbitrary Cuntz-Krieger algebras.
There have been many extensions of Cuntz-Krieger algebras, for example, Cuntz-Pimsner algebras \cite{Pi}, graph algebras and Exel-Laca algebras \cite{EL},
and there have also been many results about the ideal structures of such algebras (for example, \cite{KPW},\cite{KPRR},\cite{BPRS} and \cite{EL}).

The crossed products of $C^*$-algebras give us plenty of interesting examples and the structures of them have been examined by several authors.
In \cite{Ki}, A. Kishimoto gave a necessary and sufficient condition that the crossed products by abelian groups become simple in terms of the strong Connes spectrum.
For the case of the crossed products of the Cuntz algebras by so-called quasi-free actions of abelian groups, he gave a condition for simplicity, which is easy to check, and computed the strong Connes spectra of some of such actions.
It is hard to compute the strong Connes spectrum by its definition and there have been few examples of actions whose strong Connes spectra have been computed.

In this paper, we deal with crossed products of Cuntz algebras by quasi-free actions of arbitrary locally compact, second countable, abelian groups $G$.
The class of our algebras has many examples of simple stably projectionless $C^*$-algebras as well as AF-algebras and purely infinite $C^*$-algebras (see \cite{KK1}, \cite{KK2} or \cite{Ka}).
Our algebras may be considered as continuous counterparts of Cuntz-Krieger algebras or graph algebras.
The main purpose of this paper is to determine the ideal structures of our algebras in terms of the spectrum of the action, which is a finite subset of the dual group $\Gamma$ of $G$.
This paper is organized as follows.
After some preliminaries, we prove that the set of all ideals that are invariant under the gauge action is in a one-to-one correspondence to the set of closed subsets of the dual group $\Gamma$ of $G$ satisfying certain conditions (Theorem \ref{OneToOne}).
Next we give a necessary and sufficient condition that our algebras become simple (Theorem \ref{simple}), which gives another proof of A. Kishimoto's result.
We also give a necessary and sufficient condition that our algebras become primitive (Theorem \ref{primitive}).
In section \ref{section}, we completely determine the ideal structures of our algebras.
If actions satisfy a certain condition, which is an analogue of Condition (II) in the case of Cuntz-Krieger algebras \cite{C2} or Condition (K) in the case of graph algebras \cite{KPRR}, then one can show that all ideals are invariant under the gauge action and so one can describe all ideals in terms of closed sets of the group $\Gamma$ (Theorem \ref{idestr1}).
It is rather difficult to describe the ideal structures when actions do not satisfy the condition.
We have to determine all primitive ideals and investigate the topology of the primitive ideal spaces of our algebras.
After that, we show that when actions do not satisfy the condition, the set of all ideals corresponds bijectively to the set of closed subsets of a certain topological space satisfying a certain condition (Theorem \ref{idestr2}).
As a consequence of knowing the ideal structures completely, we can compute the strong Connes spectra of quasi-free actions on Cuntz algebras.
Our algebras can be considered as Cuntz-Pimsner algebras and by using this fact we compute the K-groups of our algebras.
Finally we conclude this paper by giving some examples and remarks.

{\bf Acknowledgment.} 
A part of this work was done while the author stayed at the Mathematical 
Sciences Research Institute, 
and he would like to express his gratitude to MSRI for their hospitality.
The author is grateful to his advisor Yasuyuki Kawahigashi for his support 
and encouragement, 
to Masaki Izumi for various comments and many suggestions.
This work was partially supported by Research Fellowship 
for Young Scientists of the Japan Society for the Promotion of Science and 
Honda Heizaemon memorial fellowship of 
the Japan Association for Mathematical Sciences.

\section{Preliminaries}\label{PRE}

In this section, we review some basic objects and fix the notation.

For $n=2,3,\ldots$, the Cuntz algebra $\On$ is the universal $C^*$-algebra 
generated by $n$ isometries $S_1,S_2,\ldots,S_n$, satisfying 
$\sum_{i=1}^nS_iS_i^*=1$ \cite{C1}.
For $k\in\N=\{0,1,\ldots\}$, we define the set $\W_n^{(k)}$ of $k$-tuples by $\W_n^{(0)}=\{\emptyset\}$ and
$$\W_n^{(k)}=\big\{ (i_1,i_2,\ldots,i_k)\ \big|\ i_j\in\{1,2,\ldots,n\}\big\}.$$
We set $\W_n=\bigcup_{k=0}^\infty\W_n^{(k)}$. 
For $\mu=(i_1,i_2,\ldots,i_k)\in\W_n$, we denote its length $k$ by $|\mu|$, and set $S_\mu=S_{i_1}S_{i_2}\cdots S_{i_k}\in\On$. 
Note that $|\emptyset|=0,\ S_\emptyset=1$. 
For $\mu=(i_1,i_2,\ldots,i_k),\nu=(j_1,j_2,\ldots,j_l)\in\W_n$, we define their product $\mu\nu\in\W_n$ by $\mu\nu=(i_1,i_2,\ldots,i_k,j_1,j_2,\ldots,j_l)$.

\medskip

We fix a locally compact abelian group $G$ which satisfies the second axiom of countability.
The dual group of $G$ is denoted by $\Gamma$ which is also a locally compact abelian group satisfying the second axiom of countability.
We always use $+$ for multiplicative operations of abelian groups except for $\T$, which is the group of the unit circle in the complex plane $\C$. 
The pairing of $t\in G$ and $\gamma\in\Gamma$ is denoted by $\ip{t}{\gamma}\in\T$. 

Let us take $\omega=(\omega_1,\omega_2,\ldots,\omega_n)\in\Gamma^n$ and fix it. 
Since the $n$ isometries $\ip{t}{\omega_1}S_1$, $\ip{t}{\omega_2}S_2$, $\ldots$, $\ip{t}{\omega_n}S_n$ also satisfy the relation above for any $t\in G$,
there is a $*$-automorphism $\alpha^\omega_t:\On\to\On$ such that $\alpha^\omega_t(S_i)=\ip{t}{\omega_i}S_i$ for $i=1,2,\ldots,n$. 
One can see that $\alpha^\omega:G\ni t\mapsto \alpha^\omega_t\in\Aut(\On)$ is a strongly continuous group homomorphism.

\bde\label{wga}
Let $\omega=(\omega_1,\omega_2,\ldots,\omega_n)\in\Gamma^n$ be given. 
We define the action $\alpha^\omega:G\curvearrowright\On$ by
$$\alpha^\omega_t(S_i)=\ip{t}{\omega_i}S_i\quad (i=1,2,\ldots,n,\ t\in G).$$
\ede

The action $\alpha^\omega:G\curvearrowright\On$ becomes quasi-free (for a definition of quasi-free actions on the Cuntz algebras, see \cite{E}).
Conversely, any quasi-free action of abelian group $G$ on $\On$ is conjugate to $\alpha^\omega$ for some $\omega\in\Gamma^n$.

By definition, the (full) crossed product $\cp$ is the universal $C^*$-algebra generated by the $*$-algebra $L^1(G,\On)$ whose multiplication and involution are defined as the following:
$$fg(t)=\int_{G} f(s)\alpha^\omega_s(g(t-s))ds\ ,\quad f^*(t)=\alpha^\omega_t(f(-t)^*),$$
for $f,g\in L^1(G,\On)$ (cf. \cite{Pe}).
The crossed product $\cp$ has a $C^*$-subalgebra $\C 1{\rtimes_{\alpha^\omega}}G$, which is isomorphic to $C_0(\Gamma)$ via the map $\C 1{\rtimes_{\alpha^\omega}}G\supset L^1(G)\ni f\mapsto\widehat{f}\in C_0(\Gamma)$, where 
$$\widehat{f}(\gamma)=\int_G \ip{t}{\gamma}f(t)dt.$$
Throughout this paper, we always consider $C_0(\Gamma)$ as a $C^*$-subalgebra of $\cp$, and use $f,g,\ldots$ for denoting elements of $C_0(\Gamma)\subset \cp$.
The Cuntz algebra $\On$ is naturally embedded into the multiplier algebra $M(\cp)$ of $\cp$. 
For each $\mu=(i_1,i_2,\ldots,i_k)$ in $\W_n$, we define an element $\omega_\mu$ of $\Gamma$ by $\omega_\mu=\sum_{j=1}^{k}\omega_{i_j}$.
For $\gamma_0\in\Gamma$, we define a (reverse) shift automorphism $\sigma_{\gamma_0}:C_0(\Gamma)\to C_0(\Gamma)$ by $(\sigma_{\gamma_0} f)(\gamma)=f(\gamma+\gamma_0)$ for $f\in C_0(\Gamma)$. 

Once noting that $\alpha^\omega_t(S_\mu)=\ip{t}{\omega_\mu}S_\mu$ for $\mu\in\W_n$, one can easily verify the following.

\bpr\label{ComRel}
For any $f\in C_0(\Gamma)\subset \cp$ and any $\mu\in\W_n$, we have $fS_\mu =S_\mu\sigma_{\omega_\mu}f$.
\epr

For a subset $X$ of a $C^*$-algebra, the linear span of $X$ is denoted by $\spa X$, and the closure of $\spa X$ is denoted by $\cspa X$.

\bpr\label{cspa}
We have $\cp=\cspa\{ S_\mu fS_\nu^*\mid \mu,\nu\in\W_n,\ f\in C_0(\Gamma)\}$.
\epr

\bprf
By Proposition \ref{ComRel}, 
$$\cspa\{ S_\mu fS_\nu^*\mid \mu,\nu\in\W_n,\ f\in C_0(\Gamma)\}=\cspa\{ S_\mu S_\nu^*f\mid \mu,\nu\in\W_n,\ f\in C_0(\Gamma)\}.$$
Obviously $\cspa\{ S_\mu S_\nu^*f\mid \mu,\nu\in\W_n,\ f\in C_0(\Gamma)\}$ contains all elements of $L^1(G,\On)$, which is dense in $\cp$. 
The proof is completed.
\eprf

We denote by $\M_k$ the $C^*$-algebra of $k \times k$ matrices for $k=1,2,\ldots$, and by $\K$ the $C^*$-algebra of compact operators of the infinite dimensional separable Hilbert space.

\section{Gauge invariant ideals}

There is an action $\beta$ of $\T$ on $\On$ called the gauge action which is defined by $\beta_t(S_i)=tS_i$ for $t\in\T, i=1,2,\ldots,n$.
We can extend this action to $\cp$ which is also called the gauge action and denoted by $\beta$.
Explicitly, $\beta_t(S_\mu fS_\nu^*)=t^{|\mu|-|\nu|}S_\mu fS_\nu^*$ for $\mu,\nu\in\W_n,\ f\in C_0(\Gamma)$ and $t\in\T$.

By an ideal, we mean a closed two-sided ideal.
In this section, we determine all the ideals which are globally invariant under the gauge action.

\bde\label{omega}
For an ideal $I$ of the crossed product $\cp$, we define the closed subset $X_I$ of $\Gamma$ by
$$X_I=\bigcap_{f\in I\cap C_0(\Gamma)}\{\gamma\in\Gamma\mid f(\gamma)=0\}.$$
\ede

In other words, $X_I$ is determined by $C_0(\Gamma\setminus X_I)=I\cap C_0(\Gamma)$ where for a closed subset $X$ of $\Gamma$, $C_0(\Gamma\setminus X)$ means the set of functions in $C_0(\Gamma)$ which vanish on $X$. 
In particular, $C_0(\Gamma\setminus \Gamma)=\{0\}$.
One can easily see that $I_1\subset I_2$ implies $X_{I_1}\supset X_{I_2}$ and $X_{I_1\cap I_2}=X_{I_1}\cup X_{I_2}$ for ideals $I_1,I_2$ of $\cp$.

\bde
A subset $X$ of $\Gamma$ is called {\em $\omega$-invariant} if $X$ is a closed set satisfying the following two conditions:
\benu
\item For any $\gamma\in X$ and any $i\in\{1,2,\ldots,n\}$, we have $\gamma+\omega_i\in X$.
\item For any $\gamma\in X$, there exists $i\in\{1,2,\ldots,n\}$ such that $\gamma-\omega_i\in X$.
\eenu
\ede

For an element $\gamma$ of an $\omega$-invariant set $X$,
one can easily show that $\gamma+\omega_\mu\in X$ for any $\mu\in\W_n$ 
and that there exists $i\in\{1,2,\ldots,n\}$ such that 
$\gamma-m\omega_{i}\in X$ for any $m\in\N$.
For a subset $X$ of $\Gamma$ and an element $\gamma_0$ of $\Gamma$, we define the subset $X+\gamma_0$ of $\Gamma$ by $X+\gamma_0=\{\gamma+\gamma_0\mid\gamma\in X\}$.
Similarly, we define $X_1+X_2=\{\gamma_1+\gamma_2\mid\gamma_1\in X_1,\ \gamma_2\in X_2\}$ for $X_1,X_2\subset\Gamma$.
A closed set $X$ is $\omega$-invariant if and only if $X=\bigcup_{i=1}^n (X+\omega_i)$.

\bpr
For any ideal $I$ of the crossed product $\cp$, the closed set $X_I$ is $\omega$-invariant.
\epr

\bprf
Take $\gamma\in X_I$ and $i\in\{1,2,\ldots,n\}$ arbitrarily.
Let $f$ be an element of $I\cap C_0(\Gamma)$.
By Proposition \ref{ComRel}, $S_{i}^*fS_{i}=S_{i}^*S_{i}\sigma_{\omega_{i}}f=\sigma_{\omega_{i}}f$.
Hence $\sigma_{\omega_{i}}f\in I\cap C_0(\Gamma)$, so we have $\sigma_{\omega_{i}}f(\gamma)=0$.
Thus, $f(\gamma+\omega_{i})=0$ for any $f\in I\cap C_0(\Gamma)$.
It implies $\gamma+\omega_{i}\in X_I$.

Let $\gamma_0$ be a point of $\Gamma$ such that $\gamma_0-\omega_i\notin X_I$ for any $i=1,2,\ldots,n$, and we will show that $\gamma_0\notin X_I$.
Since $\Gamma\setminus X_I$ is open, there is a neighborhood $U$ of $\gamma_0\in\Gamma$ such that $U-\omega_i\subset\Gamma\setminus X_I$ for any $i=1,2,\ldots,n$. 
There exists $f\in C_0(\Gamma)$ such that $f(\gamma_0)\neq 0$ and $f(\gamma)=0$ for any $\gamma\notin U$. Then $U-\omega_i\subset\Gamma\setminus X_I$ implies $\sigma_{\omega_i}f\in C_0(\Gamma\setminus X_I)\subset I$ for $i=1,2,\ldots,n$. Since
$$f=\sum_{i=1}^nS_iS_i^*f=\sum_{i=1}^nS_i\sigma_{\omega_i}fS_i^*,$$
we have $f\in I$. It implies $\gamma_0\notin X_I$.
Thus $X_I$ is $\omega$-invariant.
\eprf

We will show that for any $\omega$-invariant subset $X$, there exists a gauge invariant ideal $I$ such that $X=X_I$ (Proposition \ref{exist}).

\bde
Let $X$ be an $\omega$-invariant subset of $\Gamma$. 
We define $I_X\subset\cp$ by
$$I_X=\cspa\{S_\mu fS_\nu^*\mid \mu,\nu\in\W_n,\ f\in C_0(\Gamma\setminus X)\}.$$
\ede

\bpr
For an $\omega$-invariant subset $X$ of $\Gamma$, the set $I_X$ becomes a gauge invariant ideal of $\cp$.
\epr

\bprf
Clearly $I_X$ is a $*$-invariant closed linear space.
Since $\beta_t(S_\mu fS_\nu^*)=t^{|\mu|-|\nu|}S_\mu fS_\nu^*$ for $t\in\T$, $I_X$ is invariant under the gauge action $\beta$.
By Proposition \ref{cspa}, it suffices to show that for any $\mu_1,\nu_1,\mu_2,\nu_2\in\W_n$ and any $f\in C_0(\Gamma\setminus X),\ g\in C_0(\Gamma)$, the product $xy$ of $x=S_{\mu_1}fS_{\nu_1}^*\in I_X$ and $y=S_{\mu_2}gS_{\nu_2}^*\in\cp$ is an element of $I_X$.

If $S_{\nu_1}^*S_{\mu_2}=0$, then $xy=0\in I_X$.
Otherwise, $S_{\nu_1}^*S_{\mu_2}=S_\mu$ or $S_{\nu_1}^*S_{\mu_2}=S_\mu^*$ for some $\mu\in\W_n$.
For the case $S_{\nu_1}^*S_{\mu_2}=S_\mu$, 
$$xy=S_{\mu_1}fS_{\nu_1}^*S_{\mu_2}gS_{\nu_2}^*=S_{\mu_1}fS_\mu gS_{\nu_2}^*=S_{\mu_1\mu}(\sigma_{\omega_\mu}f)gS_{\nu_2}^*.$$
Since $f\in C_0(\Gamma\setminus X)$ and $X$ is $\omega$-invariant, we have $\sigma_{\omega_\mu}f\in C_0(\Gamma\setminus X)$. 
This implies $(\sigma_{\omega_\mu}f)g\in C_0(\Gamma\setminus X)$ and so $xy\in I_X$.
For the case $S_{\nu_1}^*S_{\mu_2}=S_\mu^*$, 
$$xy=S_{\mu_1}fS_{\nu_1}^*S_{\mu_2}gS_{\nu_2}^*=S_{\mu_1}fS_\mu^*gS_{\nu_2}^*=S_{\mu_1}f(\sigma_{\omega_\mu}g)S_{\nu_2\mu}^*.$$
Since $f\in C_0(\Gamma\setminus X)$, we have $xy\in I_X$.
It completes the proof.
\eprf

\bpr\label{exist}
For any $\omega$-invariant subset $X$ of $\Gamma$, we have $X_{I_X}=X$.
\epr

\bprf
By the definition of $I_X$, we get $X_{I_X}\subset X$.
Let us assume $X_{I_X}\subsetneqq X$.
Then there exists $f\in I_X\cap C_0(\Gamma)$ such that $f(\gamma_0)=1$ for some $\gamma_0\in X$.
Since $f\in I_X$, there exist $f_k\in C_0(\Gamma\setminus X),\ \mu_k,\nu_k\in\W_n\ (k=1,2,\ldots,K)$ such that 
$$\left\| f-\sum_{k=1}^KS_{\mu_k}f_kS_{\nu_k}^*\right\|<\frac12.$$
From this inequality, we will derive a contradiction.

Since $X$ is $\omega$-invariant, there exists $i\in\{1,2,\ldots,n\}$ such that $\gamma_0-m\omega_{i}\in X$ for any $m\in\N$.
Take $j\in\{1,2,\ldots,n\}$ with $j\neq i$. 
Set $M=\mbox{max}\{|\mu_k|,|\nu_k|\mid k=1,2,\ldots,K\}$.
Then, $S_j^*{S_i^*}^MfS_i^M S_j=\sigma_{(M\omega_i+\omega_j)}f$ and $S_j^*{S_i^*}^M S_{\mu_k}f_kS_{\nu_k}^*S_i^M S_j$ is either $0$ or $\sigma_{(m_k\omega_i+\omega_j)}f_k$ for some $m_k\leq M$.
Therefore, from 
$$\left\| S_j^*{S_i^*}^M(f-\sum_{k=1}^KS_{\mu_k}f_kS_{\nu_k}^*)S_i^M S_j\right\|<\frac12,$$
we get 
$$\left\|\sigma_{(M\omega_i+\omega_j)}f-\sum_k \sigma_{(m_k\omega_i+\omega_j)}f_k\right\|<\frac12.$$
By evaluating at $\gamma_0-M\omega_{i}-\omega_j$, we find $k\in\N$ such that $f_k(\gamma_0-(M-m_k)\omega_{i})\neq 0$. It contradicts the fact that $f_k\in C_0(\Gamma\setminus X)$ and $\gamma_0-m\omega_{i}\in X$ for any $m\in\N$.
Therefore we are done.
\eprf

By Proposition \ref{exist}, the map $I\mapsto X_I$ from the set of gauge invariant ideals $I$ of $\cp$ to the set of $\omega$-invariant subsets of $\Gamma$ is surjective.
Now, we turn to showing that this map is injective (Proposition \ref{unique}).
The method we use here is inspired by \cite{C1}.

Let $I$ be an ideal that is not $\cp$.
We investigate the quotient $(\cp)/I$ of $\cp$ by an ideal $I$. 
Since $I\cap C_0(\Gamma)=C_0(\Gamma\setminus X_I)$, a $C^*$-subalgebra $C_0(\Gamma)/(I\cap C_0(\Gamma))$ of $(\cp)/I$ is isomorphic to $C_0(X_I)$.
We will consider $C_0(X_I)$ as a $C^*$-subalgebra of $(\cp)/I$.
No confusion should occur by using same symbols $S_1,S_2,\ldots,S_n\in M((\cp)/I)$ as the ones in $M(\cp)$ for denoting the isometries of $\On$ which is naturally embedded into $M((\cp)/I)$.

For an $\omega$-invariant set $X$, we can define a $*$-homomorphism $\sigma_{\omega_\mu}:C_0(X)\to C_0(X)$ for $\mu\in\W_n$.
This map $\sigma_{\omega_\mu}$ is always surjective, 
but it is injective only in the case that $X-\omega_\mu\subset X$, 
which is equivalent to $X-\omega_\mu=X$.
One can easily verify the following.

\blem\label{cp/I}
Let $I$ be an ideal that is not $\cp$. 
\benu
\item For $\mu,\nu\in\W_n$ and $f\in C_0(X_I)$,
$S_\mu fS_\nu^*\in(\cp)/I$ is zero if and only if $f=0$.
\item For $\mu\in\W_n$ and $f\in C_0(X_I)$,
we have $fS_\mu=S_\mu \sigma_{\omega_\mu}f$.
\item $(\cp)/I=\cspa\{S_\mu fS_\nu^*\mid\mu,\nu\in\W_n,\ f\in C_0(X_I)\}$.
\eenu
\elem

We define a $C^*$-subalgebra of $(\cp)/I$, which corresponds to the AF-core for Cuntz algebras.

\bde\label{AFcore}
Let $I$ be an ideal that is not $\cp$.
We define $C^*$-subalgebras ${\mathcal F}_I^{(k)}$ ($k\in\N$) and ${\mathcal F}_I$ of $(\cp)/I$ by
\begin{align*}
{\mathcal F}_I^{(k)}&=\cspa\{S_\mu fS_\nu^*\mid \mu,\nu\in\W_n^{(k)},\ f\in C_0(X_I)\},\\
{\mathcal F}_I&=\cspa\{S_\mu fS_\nu^*\mid \mu,\nu\in\W_n,\ |\mu|=|\nu|,\ f\in C_0(X_I)\}.
\end{align*}
When $I=0$, we write simply ${\mathcal F}^{(k)},{\mathcal F}$ for ${\mathcal F}_0^{(k)},{\mathcal F}_0$.
\ede

\blem\label{F_k}
Let $I$ be an ideal that is not $\cp$. 
\benu
\item The $C^*$-subalgebra ${\mathcal F}_I^{(k)}$ of $(\cp)/I$ is isomorphic to $\M_{n^k}\otimes C_0(X_I)$ for $k\in\N$. 
\item ${\mathcal F}_I^{(k)}\subset {\mathcal F}_I^{(k+1)}$ and the inductive limit of ${\mathcal F}_I^{(k)}$ is ${\mathcal F}_I$.
\eenu
\elem

\bprf
\benu
\item Since the set $\W_n^{(k)}$ has $n^k$ elements, we may use $\{e_{\mu,\nu}\}_{\mu,\nu\in\W_n^{(k)}}$ for denoting the matrix units of $\M_{n^k}$.
For $x_1=S_{\mu_1} f_1S_{\nu_1}^*,\ x_2=S_{\mu_2} f_2S_{\nu_2}^*\in {\mathcal F}_I^{(k)}$, we have $x_1^*=S_{\nu_1}\overline{f_1}S_{\mu_1}^*$ and $x_1x_2=\delta_{\nu_1,\mu_2}S_{\mu_1} f_1f_2S_{\nu_2}^*$.
Thus the map
$$\M_{n^k}\otimes C_0(X_I)\ni e_{\mu,\nu}\otimes f \mapsto S_\mu fS_\nu^*\in {\mathcal F}_I^{(k)}$$
defines a $*$-homomorphism.
By the definition of ${\mathcal F}_I^{(k)}$, it is surjective.
By Lemma \ref{cp/I} (i), it is injective.
Thus $\M_{n^k}\otimes C_0(X_I)\cong{\mathcal F}_I^{(k)}$.
\item
Since $S_\mu fS_\nu^*=\sum_{i=1}^nS_\mu S_i \sigma_{\omega_i}fS_i^*S_\nu^*$, we have ${\mathcal F}_I^{(k)}\subset {\mathcal F}_I^{(k+1)}$.
The latter part is trivial by the definitions of ${\mathcal F}_I^{(k)}$ and ${\mathcal F}_I$.
\eprf
\eenu

\bde
A linear map $E$ from some $C^*$-algebra $A$ onto a $C^*$-subalgebra $B$ of $A$ is called a {\em conditional expectation} if $\|E\|\leq 1$ and $E(x)=x$ for any $x\in B$. 
A conditional expectation $E$ is called {\em faithful} if $E(x)=0$ implies $x=0$ for a positive element $x$ of $A$.
\ede

The following proposition essentially appeared in \cite{C1}.

\bpr\label{isom}
For $i=1,2$, let $E_i$ be a conditional expectation from a $C^*$-algebra $A_i$ onto a $C^*$-subalgebra $B_i$ of $A_i$.
Let $\varphi:A_1\to A_2$ be a $*$-homomorphism with $\varphi\circ E_1=E_2\circ\varphi$.
If the restriction of $\varphi$ on $B_1$ is injective and $E_1$ is faithful, then $\varphi$ is injective.
\epr

\bprf
Let $x$ be a positive element of $\ker\varphi\subset A_1$.
Since $\varphi\circ E_1=E_2\circ \varphi$, we have $\varphi(E_1(x))=0$.
Since $E_1(x)\in B_1$ and $\varphi$ is injective on $B_1$, we have $E_1(x)=0$.
Then $x=0$ since $E_1$ is faithful.
Thus $\ker\varphi=\{0\}$ which means that $\varphi$ is injective. 
\eprf

For an ideal $I$ which is invariant under the gauge action $\beta$, we can extend the gauge action on $\cp$ to one on $(\cp)/I$, which is also denoted by $\beta$.

\blem\label{cond.exp1}
Let $I$ be a gauge invariant ideal that is not $\cp$. 
Then, 
$$E_{I}:(\cp)/I\ni x\mapsto\int_{\T}\beta_t(x)dt\in (\cp)/I$$
is a faithful conditional expectation onto ${\mathcal F}_{I}$, where $dt$ is the normalized Haar measure on $\T$.
\elem

\bprf
For $x\in (\cp)/I$, 
$$\|E_{I}(x)\|=\left\|\int_{\T}\beta_t(x)dt\right\|\leq\int_{\T}\|\beta_t(x)\|dt=\|x\|.$$
Thus $\|E_{I}\|\leq 1$. 
One can see that $E_{I}(x)=0$ implies $x=0$ for a positive element $x\in (\cp)/I$.

For $\mu,\nu\in\W_n$ and $f\in C_0(X_I)$,
$$E_{I}(S_\mu fS_\nu^*)=\int_{\T}\beta_t(S_\mu fS_\nu^*)dt=\int_{\T}t^{|\mu|-|\nu|}(S_\mu fS_\nu^*)dt=\delta_{|\mu|,|\nu|}S_\mu fS_\nu^*.$$
Therefore $E_{I}(x)=x$ for any $x\in\spa\{S_\mu fS_\nu^*\mid \mu,\nu\in\W_n,\ |\mu|=|\nu|,\ f\in C_0(X_I)\}$, thus for any $x\in {\mathcal F}_{I}$ by the continuity of $E_{I}$.
By the above computation, $E_{I}(x)\in{\mathcal F}_{I}$ for $x\in\spa\{S_\mu fS_\nu^*\mid\mu,\nu\in\W_n,\ f\in C_0(X_I)\}$, which is dense in $(\cp)/I$ by Lemma \ref{cp/I} (iii). 
Therefore, the image of $E_{I}$ is ${\mathcal F}_{I}$ by the continuity of $E_{I}$.
We have shown that $E_{I}$ is a faithful conditional expectation onto ${\mathcal F}_{I}$.
\eprf

\bpr\label{unique}
For any gauge invariant ideal $I$, we have $I_{X_I}=I$.
\epr

\bprf
When $I=\cp$, we have $X_I=\emptyset$. Thus $I_{X_I}=\cp$.
Let $I$ be a gauge invariant ideal that is not $\cp$ and set $J=I_{X_I}$. 
By the definition, $J\subset I$. 
Hence there is a surjective $*$-homomorphism $\pi:(\cp)/J\to(\cp)/I$. 
By Lemma \ref{F_k}, the restriction of $\pi$ on ${\mathcal F}_{J}^{(k)}$ is an isomorphism from ${\mathcal F}_{J}^{(k)}$ onto ${\mathcal F}_{I}^{(k)}$ and so the restriction of $\pi$ on ${\mathcal F}_{J}$ is an isomorphism from ${\mathcal F}_{J}$ onto ${\mathcal F}_{I}$.
By Lemma \ref{cond.exp1}, there are faithful conditional expectations $E_J:(\cp)/J\to {\mathcal F}_{J}$ and $E_I:(\cp)/I\to{\mathcal F}_{I}$.
Since $E_I(\pi (x))=\pi (E_{J}(x))$ for any $x\in \spa\{S_\mu fS_\nu^*\mid \mu,\nu\in\W_n,\ f\in C_0(X)\}$, which is dense in $(\cp)/J$, we have $E_I\circ \pi=\pi\circ E_{J}$.
By Proposition \ref{isom}, $\pi$ is injective.
Therefore $I_{X_I}=I$.
\eprf

\bthm\label{OneToOne}
The maps $I\mapsto X_I$ and $X\mapsto I_X$ induce a one-to-one correspondence between the set of gauge invariant ideals of $\cp$ and the set of $\omega$-invariant subsets of $\Gamma$.
\ethm

\bprf
Combine Proposition \ref{exist} and Proposition \ref{unique}.
\eprf

\section{Simplicity and primitivity of $\cp$}

In this section, we give necessary and sufficient conditions for $\omega\in\Gamma^n$ that the crossed product $\cp$ becomes simple or primitive.

\bpr\label{cp}
Let $I$ be an ideal of $\cp$. 
Then, $I=\cp$ if and only if $X_I=\emptyset$.
\epr

\bprf
The ``only if'' part is trivial.
The ``if'' part follows from Proposition \ref{cspa}.
\eprf

For an ideal $I$ of $\cp$, we have $I_{X_I}\subset I$.
In general there exists an ideal $I$ such that $I_{X_I}\neq I$ (see Proposition \ref{bad}).
However if $X_I$ satisfies a certain condition, then $I=I_{X_I}$ (Theorem \ref{condition2}).

\bde
An $\omega$-invariant subset $X$ of $\Gamma$ is said to be {\em bad} if there exists $\gamma\in X$ such that there is only one element $i$ with $\gamma-\omega_i\in X$ in $\{1,2,\ldots,n\}$ and this element $i$ satisfies that $m\omega_i=0$ for some positive integer $m$.
An $\omega$-invariant subset $X$ of $\Gamma$ is said to be {\em good} if $X$ is not bad.
\ede

Note that $\emptyset$ is a good $\omega$-invariant set.

\blem\label{good}
An $\omega$-invariant subset $X$ of $\Gamma$ is good if and only if for any $\gamma\in X$, one of the following two conditions is satisfied:
\benu
\item There exists $i\in\{1,2,\ldots,n\}$ such that $\gamma-m\omega_i\in X$ and $\gamma-m\omega_i\neq\gamma$ for any positive integer $m$.
\item There exist $i,j\in\{1,2,\ldots,n\}$ with $i\neq j$ such that $\gamma-m\omega_i-\omega_j\in X$ for any positive integer $m$.
\eenu
\elem

\bprf
When $X$ is bad, there exists $\gamma\in X$ such that there is only one element $i$ with $\gamma-\omega_i\in X$ in $\{1,2,\ldots,n\}$ and this element $i$ satisfies that $m\omega_i=0$ for some positive integer $m$.
This $\gamma\in X$ satisfies neither condition in the statement.

Let us assume that $X$ is good and that $\gamma\in X$ does not satisfy the condition (i).
We will prove that $\gamma\in X$ satisfies the condition (ii).
Since $X$ is $\omega$-invariant, there exists $i\in\{1,2,\ldots,n\}$ such that $\gamma-m\omega_i\in X$ for any positive integer $m$.
Since $\gamma\in X$ does not satisfy the condition (i), there exists a positive integer $K$ with $K\omega_i=0$.
Since $X$ is good, there exists $j\in \{1,2,\ldots,n\}$ with $j\neq i$ such that $\gamma-\omega_j\in X$.
For any positive integer $m$, if we take $l\in\N$ so that $lK-m\geq 0$, we have
$$\gamma-\omega_j-m\omega_i=\gamma-\omega_j+(lK-m)\omega_i\in X.$$
Thus $\gamma\in X$ satisfies the condition (ii).
\eprf

\bpr\label{cond.exp2}
Let $I$ be an ideal that is not $\cp$.
If $X_I$ is a good $\omega$-invariant set, then there exists a unique conditional expectation $E_I$ from $(\cp)/I$ onto ${\mathcal F}_I$ such that $E_I(S_\mu fS_\nu^*)=\delta_{|\mu|,|\nu|}S_\mu fS_\nu^*$ for $\mu,\nu\in\W_n,\ f\in C_0(X_I)$.
\epr

\bprf
Let $\mu_l,\nu_l\in\W_n$ and $f_l\in C_0(X_I)$ be given for $l=1,2,\ldots,L$.
Then $x=\sum_{l=1}^LS_{\mu_l}f_lS_{\nu_l}^*$ is an element of $\spa\{S_\mu fS_\nu^*\mid\mu,\nu\in\W_n,\ f\in C_0(X_I)\}\subset (\cp)/I$.
Set $k=\mbox{max}\{|\mu_l|,|\nu_l|\mid l=1,2,\ldots,L\}$.
We may assume that if $|\mu_l|=|\nu_l|$, then $|\mu_l|=|\nu_l|=k$.
Let $x_0=\sum_{|\mu_l|=|\nu_l|}S_{\mu_l}f_lS_{\nu_l}^*$.
Since $x_0\in {\mathcal F}_I^{(k)}\cong\M_{n^k}\otimes C_0(X_I)$, there exists $\gamma_0\in X_I$ such that $\|x_0\|=\left\|\sum_{|\mu_l|=|\nu_l|}S_{\mu_l}f_l(\gamma_0)S_{\nu_l}^*\right\|$.
We will prove that $\|x_0\|\leq\|x\|$.

Since $X_I$ is a good $\omega$-invariant set, $\gamma_0\in X_I$ satisfies one of the two conditions in Lemma \ref{good}.
We first consider the case that $\gamma_0\in X_I$ satisfies the condition (i) in Lemma \ref{good}, that is, there exists $i\in\{1,2,\ldots,n\}$ such that $\gamma_0-m\omega_i\in X$ and $\gamma_0-m\omega_i\neq\gamma_0$ for any positive integer $m$.
We can find a neighborhood $U$ of $\gamma_0-k\omega_i\in X_I$ such that $U\cap (U+m\omega_i)=\emptyset$ for $m=1,2,\ldots,k$.
Choose a function $f$ with $0\leq f\leq 1$ satisfying $f(\gamma_0-k\omega_i)=1$ and the support of $f$ is contained in $U$.
Set $u=\sum_{\mu\in\W_n^{(k)}}S_\mu S_i^k f^{1/2} S_\mu^*\in\cp$.
Since 
$$u^*u=\sum_{\mu,\nu\in\W_n^{(k)}}S_\mu f^{1/2} {S_i^*}^k S_\mu^* S_\nu S_i^k f^{1/2}S_\nu^*=\sum_{\mu\in\W_n^{(k)}}S_\mu f S_\mu^*,$$
$u^*u$ is an element of ${\mathcal F}_I^{(k)}$ which corresponds to the element $1\otimes f$ under the isomorphism ${\mathcal F}_I^{(k)}\cong \M_{n^k}\otimes C_0(X_I)$.
Thus we have $\|u^*u\|=\sup_{\gamma\in C_0(X_I)}|f(\gamma)|=1$, and so $\|u\|=1$.
When $|\mu_l|\neq|\nu_l|$, for any $\mu,\nu\in\W_n^{(k)}$, $({S_i^*}^k S_\mu^*)S_{\mu_l} S_{\nu_l}^*(S_\nu S_i^k)$ is either zero, $S_i^m$ or ${S_i^*}^m$ with some $0<m\leq k$.
In the case that $({S_i^*}^k S_\mu^*)S_{\mu_l} S_{\nu_l}^*(S_\nu S_i^k)=S_i^m$, we have 
$$(f^{1/2} {S_i^*}^k  S_\mu^*)S_{\mu_l} S_{\nu_l}^*(S_\nu  S_i^k f^{1/2})=S_i^m (\sigma_{m\omega_1}f^{1/2})f^{1/2}=0.$$
Similarly, we have $(f^{1/2} {S_i^*}^k  S_\mu^*)S_{\mu_l} S_{\nu_l}^*(S_\nu  S_i^k f^{1/2})=0$ in the case that $({S_i^*}^k S_\mu^*)S_{\mu_l} S_{\nu_l}^*(S_\nu S_i^k)={S_i^*}^m$.
Hence if $|\mu_l|\neq|\nu_l|$, then $u^*S_{\mu_l}f_lS_{\nu_l}^*u=0$.
When $|\mu_l|=|\nu_l|=k$, we have
$$u^*S_{\mu_l} f_l S_{\nu_l}^*u=\sum_{\mu,\nu\in\W_n^{(k)}}(S_\mu f^{1/2} {S_i^*}^k S_\mu^*)S_{\mu_l} f_l S_{\nu_l}^*(S_\nu S_i^k f^{1/2} S_\nu^*)=S_{\mu_l}f (\sigma_{k\omega_i}f_l) S_{\nu_l}^*.$$
Hence $u^*xu=\sum_{|\mu_l|=|\nu_l|}S_{\mu_l}f (\sigma_{k\omega_i}f_l) S_{\nu_l}^*$.
Thus we have
$$\|u^*xu\|\geq\left\|\sum_{|\mu_l|=|\nu_l|}S_{\mu_l}\ f(\gamma_0-k\omega_i)\ \sigma_{k\omega_i}f_l(\gamma_0-k\omega_i)\ S_{\nu_l}^*\right\|=\left\|\sum_{|\mu_l|=|\nu_l|}S_{\mu_l}f_l(\gamma_0) S_{\nu_l}^*\right\|=\|x_0\|.$$
Therefore when the condition (i) is satisfied, we have $\|x_0\|\leq\|u^*xu\|\leq\|x\|$.

Next we consider the case that there exist $i,j\in\{1,2,\ldots,n\}$ with $i\neq j$ such that $\gamma_0-m\omega_i-\omega_j\in X$ for any positive integer $m$.
Set $u=\sum_{\mu\in\W_n^{(k)}}S_\mu S_i^kS_jS_\mu^*\in\On\subset M(\cp)$.
Since
$$u^*u=\sum_{\mu,\nu\in\W_n^{(k)}}(S_\mu S_j^*{S_i^*}^kS_\mu^*)(S_\nu S_i^kS_jS_\nu^*)=\sum_{\mu\in\W_n^{(k)}}S_\mu S_\mu^*=1,$$
$u$ is an isometry. Since $S_j^*{S_i^*}^kS_\mu^* S_\nu S_i^kS_j=\delta_{\mu,\nu}$ for $\mu,\nu\in\W_n$ such that $|\mu|,|\nu|\leq k$, we have $u^*S_{\mu_l} S_{\nu_l}^*u=\delta_{|\mu_l|,|\nu_l|}S_{\mu_l}S_{\nu_l}^*$, for $l=1,2,\ldots,L$.
Therefore,
\begin{align*}
u^*xu&=\sum_{l=1}^Lu^*S_{\mu_l}f_lS_{\nu_l}^*u\\
&=\sum_{l=1}^L u^*S_{\mu_l}S_{\nu_l}^*u\ (\sigma_{k\omega_i+\omega_j-\omega_{\nu_l}}f_l)\\
&=\sum_{|\mu_l|=|\nu_l|}S_{\mu_l}S_{\nu_l}^*\ (\sigma_{k\omega_i+\omega_j-\omega_{\nu_l}}f_l)\\
&=\sum_{|\mu_l|=|\nu_l|}S_{\mu_l}(\sigma_{k\omega_i+\omega_j}f_l) S_{\nu_l}^*.
\end{align*}
Since $\gamma_0-k\omega_i-\omega_j\in X_I$, we have 
$$\|u^*xu\|\geq \left\|\sum_{|\mu_l|=|\nu_l|}S_{\mu_l}\ \sigma_{k\omega_i+\omega_j}f_l(\gamma_0-k\omega_i-\omega_j)\ S_{\nu_l}^*\right\|=\left\|\sum_{|\mu_l|=|\nu_l|}S_{\mu_l}f_l(\gamma_0)S_{\nu_l}^*\right\|=\|x_0\|.$$
Therefore also for the case that the condition (ii) is satisfied, we have $\|x_0\|\leq\|x\|$.

Suppose $x$ is expressed in two ways: $x=\sum_{l=1}^LS_{\mu_l}f_lS_{\nu_l}^*=\sum_{l=1}^{L'}S_{\mu_l'}f_l'S_{\nu_l'}^*$.
Let $y=\sum_{l=1}^LS_{\mu_l}f_lS_{\nu_l}^*-\sum_{l=1}^{L'}S_{\mu_l'}f_l'S_{\nu_l'}^*$ and $y_0=\sum_{|\mu_l|=|\nu_l|}S_{\mu_l}f_lS_{\nu_l}^*-\sum_{|\mu_l'|=|\nu_l'|}S_{\mu_l'}f_l'S_{\nu_l'}^*$.
Since $\|y_0\|\leq\|y\|$ and $y=x-x=0$, we get $y_0=0$. 
Thus $\sum_{|\mu_l|=|\nu_l|}S_{\mu_l}f_lS_{\nu_l}^*=\sum_{|\mu_l'|=|\nu_l'|}S_{\mu_l'}f_l'S_{\nu_l'}^*$ which means that $x_0$ does not depend on expressions of $x$.
Hence we can define a norm-decreasing linear map $E_I$ by
\begin{align*}
E_I:\spa\{S_\mu fS_\nu^*&\mid\mu,\nu\in\W_n,\ f\in C_0(X_I)\}\ni x\\
&\mapsto\ x_0\in \spa\{S_\mu fS_\nu^*\mid\mu,\nu\in\W_n,\ |\mu|=|\nu|,\ f\in C_0(X_I)\}.
\end{align*}
Since $E_I$ is norm-decreasing and $\spa\{S_\mu fS_\nu^*\mid\mu,\nu\in\W_n,\ f\in C_0(X_I)\}$ is dense in $(\cp)/I$, we can extend $E_I$ on $(\cp)/I$ with $\|E_I\|\leq 1$ whose image is contained in ${\mathcal F}_I$.
Since $E_I(x)=x$ for $x\in \spa\{S_\mu fS_\nu^*\mid\mu,\nu\in\W_n,\ |\mu|=|\nu|,\ f\in C_0(X_I)\}$, which is dense in ${\mathcal F}_I$, we get $E_I(x)=x$ for any $x\in{\mathcal F}_I$.

Therefore $E_I$ is a conditional expectation onto ${\mathcal F}_I$.
Uniqueness follows from the condition 
$E_I(S_\mu fS_\nu^*)=\delta_{|\mu|,|\nu|}S_\mu fS_\nu^*$ 
for $\mu,\nu\in\W_n$ and $f\in C_0(X)$.
\eprf

When an ideal $I$ such that $X_I$ is good is gauge invariant, the conditional expectation $E_{I}$ defined in Proposition \ref{cond.exp2} coincides with the one in Lemma \ref{cond.exp1} by uniqueness.
Actually any ideal $I$ such that $X_I$ is good is gauge invariant.

\bthm\label{condition2}
Let $I$ be an ideal of $\cp$ such that $X_I$ is good.
Then we have $I_{X_I}=I$, and so $I$ is gauge invariant.
\ethm

\bprf
If $X_I=\emptyset$, then $I=\cp$ so $I_{X_I}=I$.
Let $I$ be an ideal that is not $\cp$, and set $J=I_{X_I}$.
By the same way as in the proof of Proposition \ref{unique}, 
there is a surjective $*$-homomorphism $\pi:(\cp)/J\to(\cp)/I$ 
whose restriction on ${\mathcal F}_{J}$ is an isomorphism from 
${\mathcal F}_{J}$ onto ${\mathcal F}_{I}$.
By Proposition \ref{cond.exp2}, there is a conditional expectation $E_I:(\cp)/I\to{\mathcal F}_{I}$.
Since $E_I(\pi (x))=\pi (E_{J}(x))$ for any $x\in \spa\{S_\mu fS_\nu^*\mid \mu,\nu\in\W_n,\ f\in C_0(X)\}$, which is dense in $(\cp)/J$, we have $E_I\circ \pi=\pi\circ E_{J}$ where $E_J:(\cp)/I\to{\mathcal F}_{J}$ is a faithful conditional expectation defined in Lemma \ref{cond.exp1}.
By Proposition \ref{isom}, $\pi$ is injective.
Therefore $I=I_{X_I}$.
\eprf

When an $\omega$-invariant set $X$ is bad, there exists an ideal $I$ with $X_I=X$ which is not gauge invariant (Proposition \ref{bad}).
As a special case of Theorem \ref{condition2}, we get the following.

\bpr\label{0}
Let $I$ be an ideal of the crossed product $\cp$. 
Then $I=0$ if and only if $X_I=\Gamma$. 
\epr

\bprf
The ``only if'' part is trivial. The ``if'' part follows from Theorem \ref{condition2} since $\Gamma$ is a good $\omega$-invariant set.
\eprf

\bde
For a non-empty subset $\I$ of $\{1,2,\ldots,n\}$, we denote by $\Omega_{\I}$ the closed semigroup generated by $\omega_1,\omega_2,\ldots,\omega_n$ and $-\omega_i$ for $i\in\I$. 
\ede

For a non-empty subset $\I$ of $\{1,2,\ldots,n\}$, the set $\Omega_{\I}$ is $\omega$-invariant.
In \cite{Ki}, A. Kishimoto found a necessary and sufficient condition for $\cp$ to become simple.
Now we can reprove it.

\bthm[Cf. {\cite[Theorem 4.4]{Ki}}]\label{simple}
The following conditions for $\omega\in\Gamma^n$ are equivalent:
\benu
\item The crossed product $\cp$ is simple.
\item Any $\omega$-invariant subset of $\Gamma$ must be $\emptyset$ or $\Gamma$.
\item $\Omega_{\{i\}}=\Gamma$ for any $i=1,2,\ldots,n$.
\eenu
\ethm

\bprf
(i)$\iff$(ii): Combine Proposition \ref{cp} and Proposition \ref{0}.

(ii)$\Rightarrow$(iii): Since $\Omega_{\{i\}}$ is a non-empty $\omega$-invariant subset, we have $\Omega_{\{i\}}=\Gamma$ for any $i$.

(iii)$\Rightarrow$(ii): Let $X$ be a non-empty $\omega$-invariant subset. 
Let us choose an element $\gamma_0\in X$. There exists $i\in\{1,2,\ldots,n\}$ such that $\gamma_0+\gamma\in X$ for any $\gamma\in\Omega_{\{i\}}$. Since $\Omega_{\{i\}}$ is $\Gamma$, we get $X=\Gamma$.
\eprf

Now we turn to determining for which $\omega\in\Gamma^n$ the crossed product $\cp$ becomes primitive.
An ideal $I$ of a $C^*$-algebra $A$ is called primitive if $I$ is a kernel of some irreducible representation.
A $C^*$-algebra $A$ is called primitive if $0$ is a primitive ideal.
When a $C^*$-algebra $A$ is separable, an ideal $I$ of $A$ is primitive if and only if $I$ is prime, i.e. $I_1\cap I_2\subset I$ implies $I_1\subset I$ or $I_2\subset I$ for ideals $I_1,I_2$ of $A$.

\bde
An $\omega$-invariant set $X$ is called {\em prime} if for any $\omega$-invariant sets $X_1,X_2$ with $X\subset X_1\cup X_2$, either $X\subset X_1$ or $X\subset X_2$ holds.
\ede

\bpr\label{prime}
If an ideal $I$ of $\cp$ is primitive, then $X_I$ becomes a prime $\omega$-invariant set.
\epr

\bprf
Let $I$ be a primitive ideal of $\cp$.
Assume that two $\omega$-invariant subsets $X_1,X_2$ of $\Gamma$ satisfy $X_I\subset X_1\cup X_2$.
Then $I_{X_1}\cap I_{X_2}\subset I_{X_I}\subset I$.
Since $I$ is prime, either $I_{X_1}\subset I$ or $I_{X_2}\subset I$.
Hence either $X_1\supset X_I$ or $X_2\supset X_I$.
Therefore $X_I$ is prime.
\eprf

In general, the converse of Proposition \ref{prime} is not true (see Corollary \ref{primitive1} and Proposition \ref{bad}).

\bpr\label{prime1}
For a non-empty $\omega$-invariant set $X$, the following are equivalent:
\benu
\item $X$ is prime.
\item For any $\gamma_0,\gamma_1\in X$ and any neighborhoods $U_0$ and $U_1$ of $\gamma_0$ and $\gamma_1$ respectively, there exist $\gamma\in X$ and $\mu,\nu\in\W_n$ such that $\gamma+\omega_{\mu}\in U_0$ and $\gamma+\omega_{\nu}\in U_1$.
\item For any $\gamma_0,\gamma_1\in X$, there exist sequences $\mu_1,\mu_2,\ldots$ and $\nu_1,\nu_2\ldots$ in $\W_n$ such that $\gamma_0-\omega_{\mu_k}, \gamma_1-\omega_{\nu_k}\in X$ for any $k$ and $\lim_{k\to\infty}\left( (\gamma_0-\omega_{\mu_k})-(\gamma_1-\omega_{\nu_k})\right)=0$.
\item $X=\gamma+\Omega_{\I}$ for some $\gamma\in\Gamma$ and non-empty $\I\subset\{1,2,\ldots,n\}$.
\eenu
\epr

\bprf
(i)$\Rightarrow$(ii): Let $X$ be a non-empty prime $\omega$-invariant set, $\gamma_0,\gamma_1$ elements of $X$, and $U_0,U_1$ neighborhoods of $\gamma_0,\gamma_1$ respectively.
Set two open sets $Y_0,Y_1$ by
\begin{align*}
Y_0&=\bigcup_{k=0}^\infty\bigcap_{\mu\in\W_n^{(k)}}\bigcup_{\nu\in\W_n}(U_0+\omega_\mu-\omega_\nu),&
Y_1&=\bigcup_{k=0}^\infty\bigcap_{\mu\in\W_n^{(k)}}\bigcup_{\nu\in\W_n}(U_1+\omega_\mu-\omega_\nu).
\end{align*}
One can easily see that $\gamma-\omega_i\in Y_0$ for any $\gamma\in Y_0$ and any $i=1,2,\ldots,n$, and that if $\gamma-\omega_i\in Y_0$ for any $i=1,2,\ldots,n$, then $\gamma\in Y_0$.
Thus the closed set $\Gamma\setminus Y_0$ is $\omega$-invariant.
Similarly, $\Gamma\setminus Y_1$ is $\omega$-invariant.
Since $\gamma_0\in Y_0$ and $\gamma_1\in Y_1$, neither $\Gamma\setminus Y_0$ nor $\Gamma\setminus Y_1$ contains $X$.
Since $X$ is prime, $(\Gamma\setminus Y_0)\cup(\Gamma\setminus Y_1)$ does not contain $X$.
Therefore we get $\gamma'\in X$ with $\gamma'\in Y_0\cap Y_1$.
Thus for $j=0,1$, there exist $k_j\in\N$ satisfying that for any $\mu_j\in\W_n^{(k_j)}$, there exists $\nu_j\in\W_n$ with $\gamma'\in U_j+\omega_{\mu_j}-\omega_{\nu_j}$.
Let $k\in\N$ be an integer with $k\geq k_0,k_1$.
Since $X$ is $\omega$-invariant, there exists $\mu'\in\W_n^{(k)}$ with $\gamma'-\omega_{\mu'}\in X$.
For $j=0,1$, there exist $\mu_j,\mu_j'\in\W_n$ with $\mu'=\mu_j\mu_j'$ and $|\mu_j|=k_j$.
Thus we get $\nu_j\in\W_n$ with $\gamma'\in U_j+\omega_{\mu_j}-\omega_{\nu_j}$ for $j=0,1$.
Set $\gamma=\gamma'-\omega_{\mu'}\in X$, $\mu=\nu_0\mu_0'$ and, $\nu=\nu_1\mu_1'$.
Then, we have $\gamma+\omega_{\mu}\in U_0$ and $\gamma+\omega_{\nu}\in U_1$.

(ii)$\Rightarrow$(iii): 
Let $\gamma_0,\gamma_1$ be elements of an $\omega$-invariant set $X$.
Let $U_1,U_2,\ldots$ be a fundamental system of neighborhoods of $0$.
From (ii), for any $k=1,2,\ldots$, there exist $\lambda_k\in X$ and $\mu_k,\nu_k\in\W_n$ such that $\lambda_k+\omega_{\mu_k}\in U_k+\gamma_0$ and $\lambda_k+\omega_{\nu_k}\in U_k+\gamma_1$.
Replacing $\{k\}$ by a subsequence if necessary, we may assume that the number of $i$'s appearing in $\mu_k$ and the one appearing in $\nu_k$ increase for any $i=1,2,\ldots,n$.
For any positive integer $k$, we have $\gamma_0-\omega_{\mu_k}\in X$ since $\gamma_0-\omega_{\mu_k}=\lim_{l\to\infty}(\lambda_l+\omega_{\mu_l}-\omega_{\mu_k})$ and $\lambda_l+\omega_{\mu_l}-\omega_{\mu_k}\in X$ when $l\geq k$.
Similarly we have $\gamma_1-\omega_{\nu_k}\in X$ for any positive integer $k$.
Since $\lim_{k\to\infty}(\gamma_0-\omega_{\mu_k}-\lambda_k)=0$ and $\lim_{k\to\infty}(\lambda_k-(\gamma_1-\omega_{\nu_k}))=0$, we have $\lim_{k\to\infty}\left( (\gamma_0-\omega_{\mu_k})-(\gamma_1-\omega_{\nu_k})\right)=0$.

(iii)$\Rightarrow$(iv): 
Take $\gamma_0\in X$ arbitrarily.
From (iii), the countable set $X'=\{\gamma\in X\mid \gamma=\gamma_0-\omega_{\mu}+\omega_\nu\mbox{ for some }\mu, \nu\in\W_n\}$ is dense in $X$.
Denote all the elements of $X'$ by $\{\lambda_1,\lambda_2,\ldots\}$.
Let $U_1,U_2,\ldots$ be a fundamental system of neighborhoods of $0$.
Let us choose a bijection $\Z^+\ni k\mapsto (m_k,l_k)\in\Z^+\times\Z^+$ where $\Z^+$ is the set of positive integers.
Thus we have $\{(\lambda_{m_k},U_{l_k})\}_{k=1}^\infty=\{(\lambda_{m},U_{l})\}_{m,l=1}^\infty$.
By (iii), for a positive integer $k$, we can recursively find $\mu_k,\nu_k\in\W_n$ satisfying that $\gamma_0-\sum_{j=1}^k\omega_{\mu_j}\in X$ and $\gamma_0-\sum_{j=1}^k\omega_{\mu_j}-(\lambda_{m_k}-\omega_{\nu_k})\in U_{l_k}$.
Since an element $\gamma_0-\sum_{j=1}^k\omega_{\mu_j}+\omega_\nu$ is in $X$ for any positive integer $k$ and any $\nu\in\W_n$, we have 
$$\overline{\left\{\left.\gamma_0-\sum_{j=1}^k\omega_{\mu_j}+\omega_\nu\ \right|\ k\in\Z^+, \nu\in\W_n \right\}}\subset X.$$
Since the set of the left hand side above contains $X'$ which is dense in $X$, the inclusion above is actually an equality.
Let $\I$ be the set of $i\in\{1,2,\ldots,n\}$ such that the number of $i$'s appearing in $\mu_1\mu_2\cdots\mu_k$ goes to infinity when $k$ goes to infinity.
For $i\notin\I$, let $m_i$ be the limit of the number of $i$'s appearing in $\mu_1\mu_2\cdots\mu_k$ when $k$ goes to infinity.
Set $\gamma=\gamma_0-\sum_{i\notin\I}m_i\omega_i$.
Then, one can see that $X=\gamma+\Omega_{\I}$.

(iv)$\Rightarrow$(i): Let $X$ be an $\omega$-invariant set such that $X=\gamma+\Omega_{\I}$ for some $\gamma\in\Gamma$ and non-empty $\I\subset\{1,2,\ldots,n\}$.
Take $\omega$-invariant sets $X_1,X_2$ with $X\subset X_1\cup X_2$.
Since $\gamma-k\left(\sum_{i\in\I}\omega_i\right)\in X$ for any positive integer $k$, either $X_1$ or $X_2$, say $X_1$, contains $\gamma-k\left(\sum_{i\in\I}\omega_i\right)$ for infinitely many $k$.
Then, $X_1$ contains $\gamma+\gamma'$ for any $\gamma'$ in the (algebraic) semigroup generated by $\omega_1,\omega_2,\ldots,\omega_n$ and $-\omega_i$ for $i\in\I$. 
Since $X_1$ is closed, $X_1\supset \gamma+\Omega_{\I}=X$.
Thus $X$ is prime.
\eprf

We will use the equivalence (i)$\iff$(iv) in Proposition \ref{prime1} most often.
The condition (ii) or (iii) in Proposition \ref{prime1} can be considered as an analogue of maximal tails in \cite{BPRS}.

\bthm\label{primitive}
The following conditions for $\omega\in\Gamma^n$ are equivalent:
\benu
\item The crossed product $\cp$ is primitive.
\item $\Gamma$ is a prime $\omega$-invariant set.
\item The closed group generated by $\omega_1,\omega_2,\ldots,\omega_n$ is equal to $\Gamma$.
\eenu
\ethm

\bprf
(i)$\Rightarrow$(ii): This follows from Proposition \ref{prime}.

(ii)$\Rightarrow$(i): It suffices to show that $0$ is prime.
Let $I_1,I_2$ be ideals of $\cp$ with $I_1\cap I_2=0$.
We have $X_{I_1}\cup X_{I_2}=X_{I_1\cap I_2}=\Gamma$.
Since $\Gamma$ is prime, either $X_{I_1}\supset\Gamma$ or $X_{I_2}\supset\Gamma$.
If $X_{I_1}\supset\Gamma$ hence $X_{I_1}=\Gamma$, then $I_1=0$ by Proposition \ref{0}.
Similarly if $X_{I_2}\supset\Gamma$, then $I_2=0$.
Thus $0$ is prime and so $\cp$ is a primitive $C^*$-algebra.

(ii)$\Rightarrow$(iii): 
By Proposition \ref{prime1}, there exist $\gamma\in\Gamma$ and non-empty 
$\I\subset\{1,2,\ldots,n\}$ with $\Gamma=\gamma+\Omega_{\I}$.
The closed group generated by $\omega_1,\omega_2,\ldots,\omega_n$ 
is equal to $\Gamma$ because it contains 
$\Omega_{\I}$ and $\Omega_{\I}=\Gamma-\gamma=\Gamma$.

(iii)$\Rightarrow$(ii): This follows from Proposition \ref{prime1} since 
$\Gamma=\Omega_{\{1,2,\ldots,n\}}$.
\eprf

One can prove the equivalence between (i) and (iii) in the above theorem
by characterization of primitivity of crossed products in terms of 
the Connes spectrum due to D. Olesen and G. K. Pedersen \cite{OP} and
the computation of the Connes spectrum of our actions $\alpha^{\omega}$ 
due to A. Kishimoto \cite{Ki}.

\section{The ideal structures of $\cp$}\label{section}

In this section, we completely determine the ideal structures of $\cp$ 
(Theorem \ref{idestr1}, Theorem \ref{idestr2}).
The ideal structures of $\cp$ depend on whether $\omega\in\Gamma^n$ satisfies 
the following condition:

\bcon\label{cond}
For each $i\in\{1,2,\ldots,n\}$, one of the following two conditions is satisfied:
\benu
\item For any positive integer $k$, $k\omega_i\neq 0$.
\item There exists $j\neq i$ such that $-\omega_j\in\Omega_{\{i\}}$.
\eenu
\econ

This condition is an analogue of Condition (II) in the case of Cuntz-Krieger algebras \cite{C2} or Condition (K) in the case of graph algebras \cite{KPRR}.
\subsection{When $\omega$ satisfies Condition \ref{cond}}

When $\omega$ satisfies Condition \ref{cond}, all ideals of $\cp$ are gauge invariant.

\bthm\label{idestr1}
When $\omega$ satisfies Condition \ref{cond}, every $\omega$-invariant set is good.
Hence any ideal is gauge invariant and there is a one-to-one correspondence between the set of ideals of $\cp$ and the set of $\omega$-invariant subsets of $\Gamma$.
\ethm

\bprf
Let $X$ be an $\omega$-invariant set and $\gamma$ be an element of $X$.
Since $X$ is $\omega$-invariant, there exists $i\in\{1,2,\ldots,n\}$ such that $\gamma+\gamma'\in X$ for any $\gamma'\in\Omega_{\{i\}}$.
If $k\omega_i\neq 0$ for any positive integer $k$, then $\gamma\in X$ satisfies the condition (i) in Lemma \ref{good} and 
if there exists $j\neq i$ such that $-\omega_j\in\Omega_{\{i\}}$, then $\gamma\in X$ satisfies the condition (ii) in Lemma \ref{good}.
Hence $X$ is a good $\omega$-invariant set.

By Theorem \ref{condition2}, any ideal $I$ of $\cp$ satisfies $I_{X_I}=I$ and is gauge invariant.
The last part follows from Theorem \ref{OneToOne}.
\eprf

\bco\label{primitive1}
When $\omega$ satisfies Condition \ref{cond}, an ideal $I$ of $\cp$ is primitive if and only if the $\omega$-invariant set $X_I$ is prime.
\eco

\bprf
It follows from Theorem \ref{idestr1}.
\eprf

\subsection{When $\omega$ does not satisfy Condition \ref{cond}}

From here until the end of this section, we assume that $\omega$ does not satisfy Condition \ref{cond}, i.e. there exists $i\in\{1,2,\ldots,n\}$ such that $k\omega_i=0$ for some positive integer $k$, and that $-\omega_j$ is not in the closed semigroup generated by $\omega_1,\omega_2,\ldots, \omega_n$ and $-\omega_i$ for any $j\neq i$.
Without loss of generality, we may assume $i=1$.
Let $K$ be the smallest positive integer satisfying $K\omega_1=0$.
Note that $-\omega_1$ is in the semigroup generated by $\omega_1,\omega_2,\ldots, \omega_n$ and that the closed set $X$ is $\omega$-invariant if and only if $X+\omega_i\subset X$ for any $i$.
Define $A_0=\spa\{S_1^kf{S_1^*}^l\mid f\in C_0(\Gamma), k,l\in\N\}$ which is a $*$-subalgebra of $\cp$ and denote its closure in $\cp$ by $A$.

\blem\label{miniproj}
For any $x\in A$, the element $(1-S_1S_1^*)x(1-S_1S_1^*)$ is of the form $(1-S_1S_1^*)f$ for some $f\in C_0(\Gamma)$.
\elem

\bprf
One can easily verify the conclusion for $x\in A_0$.
We have the conclusion for an arbitrary $x\in A$, because $\{(1-S_1S_1^*)f\mid f\in C_0(\Gamma)\}$ is closed.
\eprf

\blem\label{univ}
The $C^*$-algebra $A$ is the universal $C^*$-algebra generated by $C_0(\Gamma)$ and $S_1\in M(A)$, that is, 
for any $C^*$-algebra $B$, any $*$-homomorphism $\varphi:C_0(\Gamma)\to B$ and $u\in M(B)$ such that $u^*u=1_{M(B)}$ and $\varphi(f)u=u\varphi(\sigma_{\omega_1}f)$ for $f\in C_0(\Gamma)$, there is a unique $*$-homomorphism $\Phi:A\to B$ such that $\Phi(S_1^kf{S_1^*}^l)=u^k\varphi(f){u^*}^l$ for $k,l\in\N$ and $f\in C_0(\Gamma)$.
\elem

\bprf
Let $\widetilde{A}$ be the universal $C^*$-algebra satisfying the condition in the statement of this lemma.
We may consider $C_0(\Gamma)$ as a $C^*$-subalgebra of $\widetilde{A}$ and denote by $u\in M(\widetilde{A})$ the isometry satisfying $fu=u\sigma_{\omega_1}f$ for $f\in C_0(\Gamma)\subset \widetilde{A}$.
The $C^*$-algebra $\widetilde{A}$ is the closure of the linear span of elements $u^kf{u^*}^l$ for $k,l\in\N$ and $f\in C_0(\Gamma)$.
There is a unique $*$-homomorphism $\Psi:\widetilde{A}\to A$ such that $\Psi(u^kf{u^*}^l)=S_1^kf{S_1^*}^l$.
Since $\spa\{S_1^kf{S_1^*}^l\mid f\in C_0(\Gamma), k,l\in\N\}$ is dense in $A$, $\Psi$ is surjective.
By the universality of $\widetilde{A}$, there is an action $\widetilde{\beta}$ of $\T$ on $\widetilde{A}$ such that $\widetilde{\beta}_t(u^kf{u^*}^l)=t^{k-l}u^kf{u^*}^l$ for $t\in\T$.
Define $\widetilde{E}(x)=\int_{\T}\widetilde{\beta}_t(x)dt$ for $x\in\widetilde{A}$.
Then $\widetilde{E}$ becomes a faithful conditional expectation onto a $C^*$-subalgebra $\widetilde{B}=\cspa\{u^kf{u^*}^k\mid f\in C_0(\Gamma),k\in\N\}$ of $\widetilde{A}$.
Since $A$ is invariant under the gauge action $\beta$, we can define a conditional expectations $E$ on $A$ by $E(x)=\int_{\T}\beta_t(x)dt$.
Obviously $\Psi\circ \widetilde{E}=E\circ\Psi$.
Let us define $\widetilde{B}^{(k)}=\cspa\{u^lf{u^*}^l\mid f\in C_0(\Gamma),0\leq l\leq k\}$.
Then we have
$$\widetilde{B}^{(k)}=\left(\bigoplus_{l=0}^{k-1}\cspa\{u^l(1-uu^*)f{u^*}^l\mid f\in C_0(\Gamma)\}\right)\oplus\cspa\{u^kf{u^*}^k\mid f\in C_0(\Gamma)\}\cong\bigoplus_{l=0}^kC_0(\Gamma)$$
and $\varinjlim\widetilde{B}^{(k)}=\widetilde{B}$.
Clearly $\Psi$ is injective on $\widetilde{B}^{(k)}$, hence on $\widetilde{B}$.
By Proposition \ref{isom}, $\widetilde{A}$ is an isomorphic to $A$ via $\Psi$.
Thus $A$ is the universal $C^*$-algebra generated by $C_0(\Gamma)$ and $S_1\in M(A)$.
\eprf

\bre
The $C^*$-algebra $A$ is isomorphic to the Toeplitz algebra of the 
Hilbert module coming from the automorphism $\sigma_{\omega_1}$ of 
$C_0(\Gamma)$ \cite{Pi}, but we do not use this fact.
\ere

We will denote the elements of $\Z/K\Z$ by $0,1,\ldots,K-1$ and sometimes regard them as integers.

\bde
Let $H$ be a separable Hilbert space whose complete orthonormal system is given by $\{\xi_{k,m}\mid k\in \Z/K\Z, m\in\N \}$.
Let $p_k$ be a projection onto the subspace generated by $\{\xi_{k,m}\}_{m\in\N}$ for $k\in \Z/K\Z$ and define $u\in B(H)$ by $u(\xi_{k,m})=\xi_{k+1,m+1}$.
Let us denote by $T_K$ the $C^*$-algebra generated by $p_0,p_1,\ldots,p_{K-1}$ and $u$.
\ede

One can easily see that the elements $p_0,p_1,\ldots,p_{K-1}$ and $u$ satisfy $\sum_{k=0}^{K-1}p_k=1$, $u^*u=1$, and $p_ku=up_{k-1}$ for $k\in\Z/K\Z$,
and that $T_K=\cspa\{u^lp_k{u^*}^m\mid k\in\Z/K\Z, l,m\in\N\}$.
There is an action $\beta':\T\curvearrowright T_K$ such that $\beta'_t(u)=tu$ and $\beta'_t(p_k)=p_k$.
For $\lambda_0,\lambda_1,\ldots,\lambda_{K-1}\in\C$ and $\theta\in\T$, the diagonal matrix and the unitary 
$$\left(\begin{array}{cccc}
\lambda_0&0&\cdots&0\\
0&\lambda_1&\cdots&0\\
\vdots&\vdots&\ddots&\vdots\\
0&0&\cdots&\lambda_{K-1}
\end{array}\right),\quad
\left(\begin{array}{ccccc}0&0&\cdots&0&\theta\\
1&0&\cdots&0&0\\
0&1&\cdots&0&0\\
\vdots&\vdots&\ddots&\vdots&\vdots\\
0&0&\cdots&1&0
\end{array}\right)$$
are denoted by $\mbox{diag}\{\lambda_0,\lambda_1,\ldots,\lambda_{K-1}\}\in\M_K$ and $u_\theta\in\M_K$ respectively. 
The $C^*$-algebra $T_K$ satisfies the following.

\bpr\label{T_K}
\benu
\item For any non-zero $x\in T_K$, there exist $l,m\in\N$ with $(1-uu^*){u^*}^l x u^m(1-uu^*)\neq 0$.
\item There is a surjection $\pi:T_K\to C(\T,\M_K)$ with
$\pi\left(\sum_{k=0}^{K-1}\lambda_kp_{k}\right)(\theta)={\rm diag}\{\lambda_0,\lambda_1,\ldots,\lambda_{K-1}\}$ and $\pi(u)(\theta)=u_\theta$ for $\theta\in\T$.
\item For $t\in\T$, we define a $*$-automorphism $\beta''_t$ of $C(\T,\M_K)$ by 
$$\beta''_t(f)(\theta)={\rm diag}\{1,t,t^2,\ldots,t^{K-1}\}f(t^K\theta){\rm diag}\{1,\overline{t},\overline{t}^2,\ldots,\overline{t}^{K-1}\}$$
for $f\in C(\T,\M_K)$.
Then we have $\pi\circ\beta'_t=\beta''_t\circ\pi$.
\item If an ideal $J$ of $T_K$ satisfies that $1\in\pi(J)$, then $J=T_K$.
\eenu
\epr

\bprf
\benu
\item For any $k\in\Z/K\Z$ and $m\in\N$, $u^m(1-uu^*){u^*}^mp_k$ is the projection onto the one dimensional subspace generated by $\xi_{k,m}\in H$.
Hence, for any non-zero $x\in T_K$, there exist $k_1,k_2\in\Z/K\Z$ and $m_1,m_2\in\N$ such that $p_{k_1}u^{m_1}(1-uu^*){u^*}^{m_1} x u^{m_2}(1-uu^*){u^*}^{m_1}p_{k_2}\neq 0$.
Thus, we get $(1-uu^*){u^*}^{m_1} x u^{m_2}(1-uu^*)\neq 0$.
\item For $k\in\Z/K\Z$, set $I_k=\cspa\{u^l(1-uu^*)p_k{u^*}^m\mid l,m\in\N\}\subset T_K$.
Since $(u^l(1-uu^*)p_k{u^*}^m)(u^{l'}(1-uu^*)p_{k'}{u^*}^{m'})=\delta_{k,k'}\delta_{m,l'}u^l(1-uu^*)p_k{u^*}^{m'}$, the set $I_k$ is isomorphic to $\K$ for any $k\in\Z/K\Z$ and $I_k$ is orthogonal to $I_{k'}$ if $k\neq k'$.
One can easily see that $I=\bigoplus_{k=0}^{K-1}I_k$ becomes an ideal of $T_K$.
Let us denote by $\pi$ the quotient map from $T_K$ onto $T_K/I$.
We will prove that $T_K/I$ is isomorphic to $C(\T,\M_K)$.
Since $1-uu^*=\sum_{k=0}^{K-1}(1-uu^*)p_k\in I$, $\pi(u)$ is an unitary of $T_K/I$.
One can verify that $e_{i,j}=\pi(p_i)\pi(u)^{i-j}=\pi(u)^{i-j}\pi(p_j)$ satisfies the axiom of matrix units of $\M_K$ for $i,j\in\Z/K\Z$.
Thus $T_K/I\cong \M_K\otimes \pi(p_0)(T_K/I)\pi(p_0)$.
Since 
$$p_0 T_K p_0=\cspa\{p_0u^lp_k{u^*}^mp_0\mid k\in\Z/K\Z, l,m\in\N\}=\cspa\{p_0u^l{u^*}^mp_0\mid l,m\in\N\mbox{ with } l-m\in K\Z\},$$
$\pi(p_0)(T_K/I)\pi(p_0)=\pi(p_0 T_K p_0)$ is generated by one unitary $\pi(p_0u^Kp_0)$.
Since $p_0 T_K p_0$ and $I$ are globally invariant under the action $\beta'$ of $\T$, we can define an action $\beta'$ of $\T$ on $\pi(p_0)(T_K/I)\pi(p_0)$.
Since $\beta'_t(\pi(p_0u^Kp_0))=t^K\pi(p_0u^Kp_0)$, the spectrum of $\pi(p_0u^Kp_0)$ is $\T$.
Therefore we have 
$$\pi(p_0)(T_K/I)\pi(p_0)\cong C(\T).$$
Thus, we have $T_K/I\cong C(\T,\M_K)$ and one can easily verify that $\pi$ is a desired surjection.
\item For $k\in\Z/K\Z$, we have $\pi\circ\beta'_t(p_k)=\beta''_t\circ\pi(p_k)=\pi(p_k)$.
One can easily see that $\pi\circ\beta'_t(u)(\theta)=tu_\theta$.
On the other hand, 
\begin{align*}
\beta''_t\circ\pi(u)(\theta)&={\rm diag}\{1,t,t^2,\ldots,t^{K-1}\}\pi(u)(t^K\theta){\rm diag}\{1,\overline{t},\overline{t}^2,\ldots,\overline{t}^{K-1}\}\\
&={\rm diag}\{1,t,t^2,\ldots,t^{K-1}\}u_{t^K\theta}{\rm diag}\{1,\overline{t},\overline{t}^2,\ldots,\overline{t}^{K-1}\}\\
&=tu_\theta.
\end{align*}
Therefore we have $\pi\circ\beta'_t=\beta''_t\circ\pi$.
\item Since $1\in\pi(J)$, there exist $x_k\in I_k$ for each $k\in\Z/K\Z$ with $1-\sum_{k=0}^{K-1}x_k\in J$.
For $k\in\Z/K\Z$, there exists $y_k\in I_k$ such that $y_k\neq x_ky_k$ since $I_k$ is not unital.
For $k\in\Z/K\Z$, we have $(1-\sum_{k=0}^{K-1}x_k)y_k=y_k-x_ky_k\neq 0$ which is in $J\cap I_k$.
Since $I_k\cong\K$, $J\cap I_k\neq\{0\}$ implies $I_k\subset J$.
Thus $1\in J$ and so $J=T_K$.
\eprf
\eenu

\bpr\label{phi}
There is a unique $*$-homomorphism $\varphi:A\to C_0(\Gamma,T_K)$ such that $\varphi(S_1^lf{S_1^*}^m)=u^l\big(\sum_{k=0}^{K-1}(\sigma_{k\omega_1}f) p_k\big){u^*}^m$.
The map $\varphi$ is injective and its image is $\{f\in C_0(\Gamma,T_K)\mid f(\gamma+\omega_1)=\Phi(f(\gamma))\mbox{ for any }\gamma\in\Gamma\}\subset C_0(\Gamma,T_K)$, where $\Phi$ is a $*$-automorphism of $T_K$ satisfying $\Phi(u)=u$ and $\Phi(p_k)=p_{k-1}$.
\epr

\bprf
First note that $C_0(\Gamma)\ni f\mapsto \sum_{k=0}^{K-1}(\sigma_{k\omega_1}f) p_k\in C_0(\Gamma,T_K)$ is an injective $*$-homomorphism.
Since this map and $u\in M(C_0(\Gamma,T_K))$ satisfy the condition in Lemma \ref{univ}, there is a unique map $\varphi:A\to C_0(\Gamma,T_K)$ such that $\varphi(S_1^lf{S_1^*}^m)=u^l\big(\sum_{k=0}^{K-1}(\sigma_{k\omega_1}f) p_k\big){u^*}^m$.
As we saw in the proof of Lemma \ref{univ}, there is a faithful conditional expectation $E$ from $A$ onto the $C^*$-subalgebra $B$ of $A$ which is an inductive limit of $C_0(\Gamma)^k$.
If one defines a conditional expectation $E'$ on $C_0(\Gamma,T_K)$ by $E'(f)(\gamma)=\int_{\T}\beta'_t(f(\gamma))dt$, then one can easily see that $E,E'$ and $\varphi$ satisfy the condition in Proposition \ref{isom}.
Hence $\varphi$ is injective.
Since the $C^*$-subalgebra $\{f\in C_0(\Gamma,T_K)\mid f(\gamma+\omega_1)=\Phi(f(\gamma))\mbox{ for any }\gamma\in\Gamma\}$ is the closed linear span of $u^l\sum_{k=0}^{K-1}(\sigma_{k\omega_1}f) p_k {u^*}^m$ for $l,m\in\N$ and $f\in C_0(\Gamma)$, this subalgebra is the image of $\varphi$.
\eprf

\bde
For $\gamma\in\Gamma$, we denote by $\varphi_\gamma:A\to T_K$ the composition of the map $\varphi:A\to C_0(\Gamma,T_K)$ in Proposition \ref{phi} and the evaluation map at $\gamma\in\Gamma$.

For $(\gamma,\theta)\in\Gamma\times\T$, we define $\psi_{\gamma,\theta}:A\to\M_K$ by the composition of $\varphi_\gamma:A\to T_K$, $\pi:T_K\to C(\T,\M_K)$ in Proposition \ref{T_K} and the evaluation map at $\theta\in\T$.
\ede

Explicitly, $\varphi_\gamma(S_1^lf{S_1^*}^m)=u^l\big(\sum_{k=0}^{K-1}f(\gamma+k\omega_1)p_k\big){u^*}^m\in T_K$ and $\psi_{\gamma,\theta}(S_1^lf{S_1^*}^m)=u_\theta^l\mbox{diag}\{f(\gamma),f(\gamma+\omega_1),\ldots,f(\gamma+(K-1)\omega_1)\}{u_\theta^*}^m\in\M_K$.
As we saw in Proposition \ref{phi}, we have $\varphi_{\gamma+\omega_1}=\Phi\circ\varphi_{\gamma}$ for any $\gamma\in\Gamma$ and one can easily see that for any $(\gamma,\theta)\in\Gamma\times\T$, $\psi_{\gamma+\omega_1,\theta}(x)=u_\theta^*\psi_{\gamma,\theta}(x)u_\theta$ for $x\in A$.
For any $t\in\T$ and any $\gamma\in\Gamma$, we have $\varphi_\gamma\circ\beta_t=\beta'_t\circ\varphi_\gamma$.

Denote by $\Gamma'$ the quotient of $\Gamma$ by the subgroup generated by $\omega_1$, which is isomorphic to $\Z/K\Z$.
We denote by $[\gamma]$ and $[U]$ the images in $\Gamma'$ of $\gamma\in\Gamma$ and $U\subset\Gamma$ respectively.
We use the symbol $([\gamma],\theta)$ for denoting elements of $\Gamma'\times\T$.

\bde
For an ideal $I$ of $\cp$, we define the closed subset $Y_I$ of $\Gamma'\times\T$ by
$$Y_I=\{([\gamma],\theta)\in\Gamma'\times\T\mid \psi_{\gamma,\theta}(x)=0\mbox{ for all }x\in A\cap I\}.$$
\ede

Note that $\psi_{\gamma,\theta}(x)=0$ if and only if $\psi_{\gamma+\omega_1,\theta}(x)=0$.

\bde
A subset $Y$ of $\Gamma'\times\T$ is called {\em $\omega$-invariant} if $Y$ is a closed set satisfying that $([\gamma+\omega_i],\theta')\in Y$ for any $i\neq 1$ any $\theta'\in\T$ and any $([\gamma],\theta)\in Y$ .
\ede

To show that the closed set $Y_I$ is $\omega$-invariant for an ideal $I$, we need the following lemma.

\blem\label{lemint}
For any $x\in A$, $(\gamma,\theta)\in\Gamma\times\T$, and $i\neq 1$, we have 
$$\int_{\T}E(\psi_{\gamma+\omega_i,t}(x))dt=\lim_{m\to\infty}\psi_{\gamma,\theta}(S_i^*{S_1^*}^{mK}xS_1^{mK}S_i),$$
where $E$ is the conditional expectation from $\M_K$ to its $C^*$-subalgebra of diagonal matrices.
\elem

\bprf
Take $(\gamma,\theta)\in\Gamma\times\T$, and $i\neq 1$.
First we consider an element $x=S_1^k{S_1^*}^lf\in A_0$ for $f\in C_0(\Gamma)$ and $k,l\in\N$.
We have
$$\int_{\T}E(\psi_{\gamma+\omega_i,t}(x))dt=\int_{\T}E\left(u_t^{k-l}\psi_{\gamma+\omega_i,\theta}(f)\right)dt,$$
here note that $\psi_{\gamma+\omega_i,t}(f)$ does not depend on $t\in\T$.
When $k-l$ is not a multiple of $K$, we have $E\left(u_t^{k-l}\psi_{\gamma+\omega_i,\theta}(f)\right)=0$.
When $k-l=mK$ for some integer $m$, we have $E\left(u_t^{k-l}\psi_{\gamma+\omega_i,\theta}(f)\right)=t^m\psi_{\gamma+\omega_i,\theta}(f)$.
Since $\int_{\T}t^m\psi_{\gamma+\omega_i,\theta}(f)dt=\delta_{m,0}\psi_{\gamma+\omega_i,\theta}(f)$, we get 
$$\int_{\T}E(\psi_{\gamma+\omega_i,t}(x))dt=\delta_{k,l}\psi_{\gamma+\omega_i,\theta}(f).$$
On the other hand, for any positive integer $m$ satisfying $mK\geq k,l$, we have $S_i^*{S_1^*}^{mK}xS_1^{mK}S_i=\delta_{k,l}\sigma_{\omega_i}f\in A$, and so $\psi_{\gamma,\theta}(S_i^*{S_1^*}^{mK}xS_1^{mK}S_i)=\delta_{k,l}\psi_{\gamma+\omega_i,\theta}(f)$.
By the linearity of the equation, for any $x\in A_0$, there exists a positive integer $M$ such that 
$$\int_{\T}E(\psi_{\gamma+\omega_i,t}(x))dt=\psi_{\gamma,\theta}(S_i^*{S_1^*}^{mK}xS_1^{mK}S_i)$$
for $m\geq M$.
Approximating $x$ by elements of $A_0$, we have 
$$\int_{\T}E(\psi_{\gamma+\omega_i,t}(x))dt=\lim_{m\to\infty}\psi_{\gamma,\theta}(S_i^*{S_1^*}^{mK}xS_1^{mK}S_i)$$
for an arbitrary $x\in A$.
\eprf

\bpr
For any ideal $I$ of $\cp$, the closed set $Y_I$ is $\omega$-invariant.
\epr

\bprf
Take $([\gamma],\theta)\in Y_I$, $i\neq 1$ and $\theta'\in\T$.
By Lemma \ref{lemint}, for any positive element $x$ of $A\cap I$, 
$$\int_{\T}E(\psi_{\gamma+\omega_i,t}(x))dt=\lim_{m\to\infty}\psi_{\gamma,\theta}(S_i^*{S_1^*}^{mK}xS_1^{mK}S_i)=0$$
since $S_i^*{S_1^*}^{mK}xS_1^{mK}S_i\in A\cap I$.
Hence $E(\psi_{\gamma+\omega_i,\theta'}(x))=0$.
Since $E$ is faithful, $\psi_{\gamma+\omega_i,\theta'}(x)=0$ for any $x\in A\cap I$.
It implies that $([\gamma+\omega_i],\theta')\in Y_I$.
Thus $Y_I$ is $\omega$-invariant.
\eprf

For an ideal $I$ of $\cp$, the closed set $X_I$, defined in Definition \ref{omega}, is determined from $Y_I$ as follows:

\bpr\label{YtoX}
For an ideal $I$ of $\cp$, 
we have $X_I=\{\gamma\in\Gamma\mid ([\gamma],\theta)\in Y_I\mbox{ for some } \theta\in\T\}$.
\epr

\bprf
When $\gamma\notin X_I$, there exists $f\in C_0(\Gamma)\cap I\subset A\cap I$ with $f(\gamma)\neq 0$.
Then for any $\theta\in\T$, $\psi_{\gamma,\theta}(f)\neq 0$.
Thus $([\gamma],\theta)\notin Y_I$ for any $\theta\in\T$.
Conversely assume $\gamma\in\Gamma$ satisfies $([\gamma],\theta)\notin Y_I$ for any $\theta\in\T$.
Then the ideal $J=\varphi_\gamma(A\cap I)$ of $T_K$ satisfies $1\in \pi(J)$ where $\pi$ is the surjection in Proposition \ref{T_K}.
By Proposition \ref{T_K} (iv), we have $\varphi_\gamma(A\cap I)=T_K$.
Hence there exists $x\in A\cap I$ with $\varphi_\gamma(x)=1$.
Since $A$ is isomorphic to the subalgebra of $C_0(\Gamma,T_K)$, there exists $y\in A$ such that $xy\in C_0(\Gamma)$ and $\varphi_{\gamma}(xy)=1$.
Since $xy\in C_0(\Gamma)\cap I$, we have $\gamma\notin X_I$.
\eprf

We get the $\omega$-invariant set $Y_I$ from an ideal $I$ of $\cp$.
Conversely, from an $\omega$-invariant set $Y$, we can construct the ideal $I_Y$ of $\cp$.

\bde
Let $Y$ be an $\omega$-invariant subset of $\Gamma'\times\T$. 
We define the subset $Y|_{\Gamma}$ of $\Gamma$ by
$$Y|_{\Gamma}=\{\gamma\in\Gamma\mid\mbox{there exist }i\neq 1\mbox{ and }\theta\in\T\mbox{ such that }([\gamma-\omega_i],\theta)\in Y\},$$
\ede

Since $\T$ is compact, the set $X_i=\{\gamma\in\Gamma\mid \mbox{there exists }\theta\in\T\mbox{ such that }([\gamma-\omega_i],\theta)\in Y\}$ is closed for $i=2,3,\dots, n$.
Thus $Y|_{\Gamma}=\bigcup_{i=2}^n X_i$ is a closed set of $\Gamma$.
Since $Y$ is $\omega$-invariant, we have $([\gamma],\theta)\in Y$ for any $\gamma\in Y|_{\Gamma}$ and any $\theta\in\T$.

\bde
For an $\omega$-invariant subset $Y$ of $\Gamma'\times\T$, 
we define $J_Y\subset A$ and $I_Y\subset\cp$ by
$$J_Y=\{x\in A\mid \psi_{\gamma,\theta}(x)=0\mbox{ for }([\gamma],\theta)\in Y,\mbox{ and } \varphi_\gamma(x)=0\mbox{ for }\gamma\in Y|_{\Gamma}\},$$
$$I_Y=\cspa\{S_\mu xS_\nu^*\mid \mu,\nu\in\W_n,\ x\in J_Y\}.$$
\ede

Clearly by the definition, $J_Y$ is an ideal of $A$.
To see that $I_Y$ is an ideal of $\cp$, we need the following lemmas.

\blem\label{J_Y1}
For any $x\in A$ and $i\neq 1$, we have $xS_i=\sum_{k=0}^{\infty}S_1^kS_i(S_i^*({S_1^*}^kx)S_i)$.
\elem

\bprf
Let us take $x\in A$ and $i\neq 1$.
When $j$ is not $1$ nor $i$, we have $S_j^*xS_i=0$ because $S_j^*x_0S_i=0$ for any element $x_0$ of $A_0$.
From this, one can show by the induction with respect to $m$ that 
$$xS_i-\sum_{k=0}^{m-1}S_1^kS_i(S_i^*({S_1^*}^kx)S_i)=S_1^m{S_1^*}^mxS_i$$
 for all positive integer $m$.
We have $\lim_{m\to\infty}S_1^m{S_1^*}^mxS_i=0$,
since for any element $x_0$ of $A_0$, we have $S_1^m{S_1^*}^mx_0S_i=0$ if we choose $m$ sufficiently large. 
It completes the proof.
\eprf

\blem\label{J_Y2}
Let $Y$ be an $\omega$-invariant subset of $\Gamma'\times\T$. 
For any $x\in J_Y$ and $i\neq 1$, we have $S_i^* xS_i\in J_Y$.
\elem

\bprf
Since $S_i=(1-S_1S_1^*)S_i$, we have $S_i^* xS_i=S_i^* (1-S_1S_1^*)x(1-S_1S_1^*)S_i$.
By Lemma \ref{miniproj}, $(1-S_1S_1^*)x(1-S_1S_1^*)=(1-S_1S_1^*)f$ for some $f\in C_0(\Gamma)$.
Hence $S_i^* xS_i=\sigma_{\omega_i}f$.
Since $x$ is in $J_Y$, so is $(1-S_1S_1^*)f$.
Let $(\gamma,\theta)\in Y$ be given.
Since $\gamma+\omega_i\in Y|_{\Gamma}$, we have $\varphi_{\gamma+\omega_i}((1-S_1S_1^*)f)=0$.
It implies that $(1-uu^*)\sum_{k=0}^{K-1}f(\gamma+\omega_i+k\omega_1)p_k=0$, and so $f(\gamma+\omega_i+k\omega_1)=0$ for $k=0,1,\ldots,K-1$.
Therefore, we have
$$\psi_{\gamma,\theta}(S_i^* xS_i)=\psi_{\gamma,\theta}(\sigma_{\omega_i}f)=\mbox{diag}\{f(\gamma+\omega_i),f(\gamma+\omega_i+\omega_1),\ldots,f(\gamma+\omega_i+(K-1)\omega_1)\}=0.$$
Similarly, if $\gamma\in Y|_{\Gamma}$, then $\gamma+\omega_i\in Y|_{\Gamma}$, and so
$$\varphi_\gamma(S_i^* xS_i)=\varphi_\gamma(\sigma_{\omega_i}f)=\sum_{k=0}^{K-1}f(\gamma+\omega_i+k\omega_1)p_k=0.$$
Therefore $S_i^* xS_i\in J_Y$.
\eprf

\bpr
For an $\omega$-invariant subset $Y$ of $\Gamma'\times\T$, 
$I_Y$ is an ideal of $\cp$.
\epr

\bprf
By definition, $I_Y$ is a $*$-invariant closed subspace of $\cp$.
To prove that $I_Y$ is an ideal of $\cp$, it suffices to show that for any $\mu,\nu\in\W_n$, $x\in J_Y$, the products of $y=S_\mu xS_\nu^*\in I_Y$ and $S_i$, $S_i^*$ ($i=1,2,\ldots,n$), or $f\in C_0(\Gamma)$ are in $I_Y$.
It is clear that $yS_i^*=S_\mu xS_\nu^*S_i^*\in I_Y$ and $yf=S_\mu x(\sigma_{\omega_\nu}f) S_\nu^*\in I_Y$.
It is also clear that $yS_i\in I_Y$ when $\nu\neq\emptyset$ or $i=1$.
Hence all we have to do is to prove $S_\mu xS_i\in I_Y$ for $\mu\in\W_n$, $x\in J_Y$ and $i\neq 1$.
By Lemma \ref{J_Y1}, we have $S_\mu xS_i=\sum_{k=0}^{\infty}S_\mu S_1^kS_i(S_i^*({S_1^*}^kx)S_i)$.
By Lemma \ref{J_Y1}, we have $S_i^*({S_1^*}^kx)S_i\in J_Y$ for any positive integer $k$.
Therefore $S_\mu xS_i\in I_Y$.
We are done.
\eprf

We will show that an ideal $I_Y$ satisfies that $Y_{I_Y}=Y$ for any $\omega$-invariant subset $Y$ of $\Gamma'\times\T$.

\blem\label{Y_I_Y1}
For an $\omega$-invariant subset $Y$ of $\Gamma'\times\T$, we have $I_Y\cap A=J_Y$.
\elem

\bprf
By the definition of $I_Y$, we have $I_Y\cap A\supset J_Y$.
We will prove the other inclusion.
Take $x\in I_Y\cap A$.
For an arbitrary $\e>0$, there exist $\mu_l,\nu_l\in\W_n$ and $x_l\in J_Y$ for $l=1,2,\ldots,L$ such that $\left\|x-\sum_{l=1}^LS_{\mu_l}x_lS_{\nu_l}^*\right\|<\e$.
Take a positive integer $m$ such that $m\geq |\mu_l|,|\nu_l|$ for $l=1,2,\ldots,L$.
Then, $\left\|{S_1^*}^mxS_1^m-\sum_{l=1}^Lx_l'\right\|<\e$ where $x_l'={S_1^*}^mS_{\mu_l}x_lS_{\nu_l}^*S_1^m$ for $l=1,2,\ldots,L$.
Since $x_l'\in J_Y$, we have $\|\psi_{\gamma,\theta}({S_1^*}^mxS_1^m)\|<\e$ for $([\gamma],\theta)\in Y$.
Since $\psi_{\gamma,\theta}(S_1)$ is a unitary, we have $\|\psi_{\gamma,\theta}(x)\|<\e$ for arbitrary $\e>0$.
Hence, we have $\psi_{\gamma,\theta}(x)=0$ for any $([\gamma],\theta)\in Y$.

Let $\gamma$ be an element of $\Gamma$.
Assume $\varphi_\gamma(x)\neq 0$ and we will prove that $\gamma\notin Y|_{\Gamma}$.
By Proposition \ref{T_K} (i), there exist $k,l\in\N$ satisfying $(1-uu^*){u^*}^k\varphi_\gamma(x)u^l(1-uu^*)\neq 0$.
Set $y=(1-S_1S_1^*){S_1^*}^kxS_1^l(1-S_1S_1^*)\in I_Y$.
Then there exists $f\in C_0(\Gamma)$ with $y=(1-S_1S_1^*)f$.
Since $\varphi_\gamma(y)\neq 0$, there exists an integer $k$ with $0\leq k\leq K-1$ such that $f(\gamma+k\omega_1)\neq 0$.
Therefore, for any $i\neq 1$ and any $\theta\in\T$, we have $\psi_{\gamma-\omega_i,\theta}(S_i^*yS_i)=\psi_{\gamma-\omega_i,\theta}(\sigma_{\omega_i}f)\neq 0$.
By the former part of this proof, we have $([\gamma-\omega_i],\theta)\notin Y$ for any $i\neq 1$ and any $\theta\in\T$ because $S_i^*yS_i\in I_Y\cap A$.
It implies that $\gamma\notin Y|_{\Gamma}$.

Therefore $x\in J_Y$.
We have proved the inclusion $I_Y\cap A\subset J_Y$ and so $I_Y\cap A=J_Y$.
\eprf

\blem\label{Y_I_Y2}
Let $Y$ be an $\omega$-invariant subset of $\Gamma'\times\T$. 
For any $([\gamma_0],\theta_0)\notin Y$, there exists $x_0\in J_Y$ such that $\psi_{\gamma_0,\theta_0}(x_0)\neq 0$.
\elem

\bprf
Since $([\gamma_0],\theta_0)\notin Y$, we have $\gamma_0\notin Y|_{\Gamma}$.
Since $Y|_{\Gamma}$ and $Y$ are closed, there exist a neighborhood $U\subset\Gamma$ of $\gamma_0$ and a neighborhood $V\subset\T$ of $\theta_0$ such that $U+\omega_1=U$, $Y|_{\Gamma}\cap U=\emptyset$ and $Y\cap ([U]\times V)=\emptyset$.
Take $g\in C(\T)$ whose support is contained in $V$ and satisfying $g(\theta_0)=1$.
The $C^*$-algebra $C^*(S_1^K)$ generated by $S_1^K$ in $M(A)$ is isomorphic to the Toeplitz algebra.
There exists an element $x\in C^*(S_1^K)\subset M(A)$ such that $\Psi_\theta(x)=g(\theta)$ where $\Psi_\theta$ is the $*$-homomorphism from $C^*(S_1^K)$ to $\C$ determined by $\Psi_\theta(S_1^K)=\theta$ for $\theta\in\T$.
Take $f\in C_0(\Gamma)$ whose support is contained in $U$ and satisfying $f(\gamma_0)=1$ and set $x_0=xf\in A$.
We have $x_0\in J_Y$ since $\varphi_{\gamma}(x_0)=0$ when $\gamma\notin U$ and $\psi_{\gamma,\theta}(x_0)=0$ when $(\gamma,\theta)\notin U\times V$.
Thus we get $x_0\in J_Y$ with $\psi_{\gamma_0,\theta_0}(x_0)\neq 0$.
\eprf

\bpr\label{exist2}
For any $\omega$-invariant subset $Y$ of $\Gamma'\times\T$, we have $Y_{I_Y}=Y$.
\epr

\bprf
Combine Lemma \ref{Y_I_Y1} and Lemma \ref{Y_I_Y2}.
\eprf

By Proposition \ref{exist2}, the map $I\mapsto Y_I$ from the set of ideals $I$ of $\cp$ to the set of $\omega$-invariant subsets of $\Gamma'\times\T$ is surjective.
We will prove this map is injective in the next subsection.
We conclude this subsection by proving some results on $I_Y$.

\bpr\label{rotation}
Let $Y$ be an $\omega$-invariant subset of $\Gamma'\times\T$.
For $t\in\T$, set $\rho_t(Y)=\{([\gamma],\theta)\in\Gamma'\times\T\mid ([\gamma],t\theta)\in Y\}$.
Then $\rho_t(Y)$ becomes $\omega$-invariant and $\beta_t(I_Y)=I_{\rho_{t^K}(Y)}$ where $\beta$ is the gauge action.
In particular, if $t^K=1$ then $\beta_t(I_Y)=I_Y$.
\epr

\bprf
By Proposition \ref{T_K} (iii), $\psi_{\gamma,\theta}(\beta_t(x))=0$ if and only if $\psi_{\gamma,t^K\theta}(x)=0$ for $x\in A$.
Since $\rho_{t^K}(Y)|_{\Gamma}=Y|_{\Gamma}$, we have
\begin{align*}
\beta_t(J_Y)&=\{x\in A\mid \psi_{\gamma,\theta}(\beta_t(x))=0\mbox{ for }([\gamma],\theta)\in Y, \varphi_\gamma(\beta_t(x))=0\mbox{ for }\gamma\in Y|_{\Gamma}\}\\
&=\{x\in A\mid \psi_{\gamma,t^K\theta}(x)=0\mbox{ for }([\gamma],\theta)\in Y, \varphi_\gamma(\beta_t(x))=0\mbox{ for }\gamma\in\rho_{t^K}(Y)|_{\Gamma}\}\\
&=J_{\rho_{t^K}(Y)}.
\end{align*}
Hence, $\beta_t(I_Y)=I_{\rho_{t^K}(Y)}$.
\eprf

A relation between $I_Y$ and $I_X$ is the following.

\bpr\label{I_YandI_X}
\benu
\item For an $\omega$-invariant subset $X$ of $\Gamma$,
$Y=[X]\times\T$ is an $\omega$-invariant subset of $\Gamma'\times\T$ and $I_Y=I_X$.
\item For an $\omega$-invariant subset $Y$ of $\Gamma'\times\T$,
$X=\{\gamma\in\Gamma\mid ([\gamma],\theta)\in Y \mbox{ for some }\theta\in\T\}$ is an $\omega$-invariant subset of $\Gamma$ and $I_{X}=\bigcap_{t\in\T}\beta_t(I_Y)$.
\eenu
\epr

\bprf
\benu
\item It is easy to see that $Y=[X]\times\T$ becomes an $\omega$-invariant set.
Noting that $Y_{I_Y}=Y$ by Proposition \ref{exist2}, we have $X_{I_Y}=X$ from Proposition \ref{YtoX}.
By Proposition \ref{rotation}, $I_Y$ is a gauge invariant ideal of $\cp$.
Therefore $I_Y=I_X$.
\item 
Noting that $Y_{I_Y}=Y$, we have $X_{I_Y}=X$ from Lemma \ref{Y_I_Y1}.
Hence $\big(\bigcap_{t\in\T}\beta_t(I_Y)\big)\cap C_0(\Gamma)=C_0(\Gamma\setminus X)=\bigcap_{t\in\T}\beta_t(I_Y\cap C_0(\Gamma))=C_0(\Gamma\setminus X)$.
Since the ideal $\bigcap_{t\in\T}\beta_t(I_Y)$ is gauge invariant,
we have $I_{X}=\bigcap_{t\in\T}\beta_t(I_Y)$.
\eprf
\eenu

\bpr\label{bad}
Let $X$ be an $\omega$-invariant subset of $\Gamma$ and set $X'=\bigcup_{i=2}^n (X+\omega_i)$ which becomes an $\omega$-invariant set satisfying $X'\subset X$.
The set $X$ is a bad $\omega$-invariant set if and only if $X'\subsetneqq X$.
When $X$ is bad, the set $Y=([X]\times\{1\})\cup ([X']\times\T)$ becomes an $\omega$-invariant subset of $\Gamma'\times\T$ satisfying $Y\subsetneqq [X]\times\T$.
Any closed set $Y'$ satisfying $Y\subset Y'\subset [X]\times\T$ is $\omega$-invariant and satisfies $X_{I_{Y'}}=X$.
\epr

\bprf
If $X$ is good, then for any $\gamma\in X$ there exists $i\neq 1$ with $\gamma-\omega_i\in X$.
Hence $X'=X$.
Conversely, if $X'=X$, then for any $\gamma\in X$ there exists $i\neq 1$ with $\gamma-\omega_i\in X$.
Hence $\gamma\in X$ satisfies the condition (ii) in Lemma \ref{good}.
Therefore $X$ is good.
When $X'\subsetneqq X$, it is easy to see that any closed set $Y'$ satisfying $Y\subset Y'\subset [X]\times\T$ is $\omega$-invariant.
The last statement follows from Proposition \ref{YtoX}.
\eprf

By Proposition \ref{bad}, we can find many ideals $I$ with $X_I=X$ if $X$ is a bad $\omega$-invariant subset of $\Gamma$ (note that a bad $\omega$-invariant set exists whenever $\omega$ does not satisfy Condition \ref{cond}).

\subsection{Primitive ideals}

Now, we turn to showing that $I_{Y_I}=I$ for any ideal $I$ (Theorem \ref{idestr2}).
To see this, we examine the primitive ideal space of $\cp$.
Let $P$ be a primitive ideal of $\cp$.
By Proposition \ref{prime}, the closed set $X_P$ of $\Gamma$ is prime. 
Hence there exist non-empty subset $\I$ of $\{1,2,\ldots,n\}$ and $\gamma_0\in\Gamma$ such that $X_P=\gamma_0+\Omega_{\I}$ by Proposition \ref{prime1}.

\bpr
If a non-empty subset $\I$ is not $\{1\}$, 
then $X=\gamma_0+\Omega_{\I}$ is a good $\omega$-invariant set for any $\gamma_0\in\Gamma$.
\epr

\bprf
There is $i\in\I$ which is not $1$.
For any $\gamma\in X$, we have $\gamma-m\omega_1-\omega_i\in X$ for any positive integer $m$.
Hence, $X$ is good by Lemma \ref{good}.
\eprf

If a primitive ideal $P$ satisfies $X_P=\gamma_0+\Omega_{\I}$ for $\I\neq\{1\}$, then $P=I_{X_P}$ by Theorem \ref{condition2}.
Conversely for any $\gamma_0\in\Gamma$ and $\I\neq\{1\}$, the ideal $I_{\gamma_0+\Omega_{\I}}$ is primitive.

\blem\label{primideal1}
Let $\gamma_0\in\Gamma$, $\I\neq\{1\}$, and $X=\gamma_0+\Omega_{\I}$.
If an ideal $I$ satisfies $X\subset X_I$, then $I\subset I_X$.
\elem

\bprf
Let us write $P=I_X$.
We will first show that $X_{I+P}=X$.
Clearly, $X_{I+P}\subset X$.
To derive a contradiction, let us assume $X_{I+P}\neq X$.
Choose $\gamma\in X$ with $\gamma\notin X_{I+P}$.
Then there exists $f\in (I+P)\cap C_0(\Gamma)$ with $f\geq 0$ and $f(\gamma)=1$.
Let us denote $f=x_1+y_1$ where $x_1\in I$ and $y_1\in P$.
Take $i\in\I$ with $i\neq 1$.
For $m\in\N$, let us define $u_m=\sum_{\mu\in\W_n^{(m)}}S_\mu S_1^{mK}S_iS_\mu^*\in M(\cp)$.
We have $u_m^*u_m=1$ for any $m\in\N$.
By checking for elements which are finite sums of $S_\mu fS_\nu^*$,
one can prove $\lim_{m\to\infty}u_m^* x u_m=\sigma_{\omega_i}(E(x))$ for any $x\in\cp$ where $E$ is the conditional expectation onto ${\mathcal F}$ and $\sigma_{\omega_i}$ is an automorphism of ${\mathcal F}$ coming from the shift of ${\mathcal F}^{(k)}\cong C_0(\Gamma)\otimes\M_{n^k}$.
Set $x_2=\sigma_{\omega_i}(E(x_1))$ and $y_2=\sigma_{\omega_i}(E(y_1))$.
Then we have $x_2\in I$, $y_2\in P$ and $\sigma_{\omega_i}(f)=x_2+y_2$.
For sufficiently large integer $k$ which is a multiple of $K$, one can find $x_3\in I\cap {\mathcal F}^{(k)}, y_3\in P\cap{\mathcal F}^{(k)}$ with $\|x_3-x_2\|<1/2, \|y_3-y_2\|<1/2$.
Note that $I\cap {\mathcal F}^{(k)}\cong C_0(\Gamma\setminus X_I)\otimes\M_{n^{k}}$ and $P\cap {\mathcal F}^{(k)}\cong C_0(\Gamma\setminus X)\otimes\M_{n^{k}}$.
Let $\varphi$ be an evaluation map at $\gamma-\omega_i\in\Gamma$ from ${\mathcal F}^{(k)}\cong C_0(\Gamma)\otimes\M_{n^{k}}$ to $\M_{n^{k}}$.
Since $\gamma-\omega_i\in X\subset X_I$, we have $\varphi(x_3)=0$ and $\varphi(y_3)=0$.
Since $\sigma_{\omega_i}(f)=\sum_{\mu\in\W_n^{(k)}}S_{\mu}\sigma_{\omega_i+\omega_\mu}(f)S_{\mu}^*\in {\mathcal F}^{(k)}$, we have $\varphi(\sigma_{\omega_i}(f))=\sum_{\mu\in\W_n^{(k)}}f(\gamma+\omega_\mu)S_{\mu}S_{\mu}^*$.
Since $f\geq 0$, we have $\varphi(\sigma_{\omega_i}(f))\geq f(\gamma+k\omega_{1})S_1^k{S_1^*}^k=S_1^k{S_1^*}^k$.
Therefore $\|\varphi(\sigma_{\omega_i}(f))\|\geq 1$ which contradicts the fact that $\|\sigma_{\omega_i}(f)-x_3-y_3\|<1$.
Hence $X_{I+P}=X$.

Since $X$ is a good $\omega$-invariant set, $X_{I+P}=X$ implies $I+P=I_X=P$.
Therefore, $I\subset P$.
\eprf

\bpr
Let $\gamma_0\in\Gamma$, $\I\neq\{1\}$, and $X=\gamma_0+\Omega_{\I}$.
The ideal $P=I_X$ is primitive.
\epr

\bprf
Since $\cp$ is separable, it suffices to show that $P$ is prime.
Let $I_1,I_2$ be ideals of $\cp$ with $I_1\cap I_2\subset P$.
Then $X_{I_1}\cup X_{I_2}\supset X$.
Since $X$ is prime, either $X_{I_1}\supset X$ or $X_{I_2}\supset X$.
By Lemma \ref{primideal1}, we have either $I_1\subset P$ or $I_2\subset P$.
Thus, $P$ is prime.
\eprf

Next we will determine all primitive ideals $P$ with $X_P=\gamma_0+\Omega_{\{1\}}$ for some $\gamma_0\in\Gamma$.
Note that $\gamma_0+\Omega_{\{1\}}$ is a bad $\omega$-invariant set.
Let $\Omega=\{\omega_{\mu}\mid \mu\in\W_n\}$ which is the semigroup generated by $\omega_1,\omega_2,\ldots,\omega_n$.
Using that $\omega\in\Gamma^n$ does not satisfy Condition \ref{cond}, we can show that $\Omega$ has no accumulation point.

\bpr\label{disc}
For any $\gamma\in\Gamma$, there exists a neighborhood $U$ of $\gamma$ with $(U\setminus\{\gamma\})\cap\Omega=\emptyset$.
\epr

\bprf
To derive a contradiction, assume that there exists $\gamma\in\Gamma$ such that $(U\setminus\{\gamma\})\cap\Omega\neq\emptyset$ for any neighborhood $U$ of $\gamma$.
One can find $\mu_0,\mu_1,\ldots$, $\mu_k$, $\ldots\in\W_n$ with $\lim_{k\to\infty}\omega_{\mu_k}=\gamma$ and $\omega_{\mu_k}\neq\gamma$ for any $k\in\N$.
Replacing $\{k\}$ by a subsequence if necessary, we may assume the number of $i$'s appearing in $\mu_k$ does not decrease for any $i=2,3,\ldots,n$. 
There exists $i_0\neq 1$ such that the number of $i_0$'s appearing in $\mu_k$ goes to infinity since $\omega_{\mu_k}\neq\gamma$ for any $k\in\N$.
Replacing $\{k\}$ by a subsequence if necessary, we may assume the number of $i_0$'s appearing in $\mu_k$ increases strictly.
We get $\lim_{k\to\infty}(\omega_{\mu_k}-\omega_{\mu_{k-1}}-\omega_{i_0})=\gamma-\gamma-\omega_{i_0}=-\omega_{i_0}$.
Since $\omega_{\mu_k}-\omega_{\mu_{k-1}}-\omega_{i_0}\in\Omega\subset\Omega_{\{1\}}$, we have $-\omega_{i_0}\in\Omega_{\{1\}}$.
It contradicts the assumption for $\omega$.
\eprf

\bco\label{disc2}
We have $\Omega_{\{1\}}=\Omega$ and $\Omega$ is a discrete set.
\eco

By the corollary above, we can define the following.

\bde\label{DefP}
For $([\gamma],\theta)\in\Gamma'\times\T$, we set $Y_{([\gamma],\theta)}=\big\{([\gamma],\theta)\big\}\cup\big(\left([\gamma+\Omega]\setminus\{[\gamma]\}\right)\times\T\big)$ which is an $\omega$-invariant closed subset of $\Gamma'\times\T$.
We write $P_{([\gamma],\theta)}$ for denoting $I_{Y_{([\gamma],\theta)}}$.
\ede

We will show that $P_{([\gamma],\theta)}$ is a primitive ideal for any $([\gamma],\theta)\in\Gamma'\times\T$.
To see this, we need the following proposition.
This will be used to determine the topology of primitive ideal space of $\cp$.
Let us define a subset $\W_n^+$ of $\W_n$ by $\W_n^+=\{(i_1,i_2,\ldots,i_k)\in\W_n\mid i_k\neq 1\}\cup\{\emptyset\}$. 

\bpr\label{local}
Let $X$ be a compact subset of $\Gamma$ such that $X\cap (X+\gamma)=\emptyset$ for any $\gamma\in\Omega\setminus\{0\}$.
If we set $X_1=X+\Omega$ and $X_2=X+(\Omega\setminus\{0,\omega_1,\ldots,(K-1)\omega_1\})$, then $X_1$ and $X_2$ become $\omega$-invariant sets and $I_{X_2}/I_{X_1}\cong\K\otimes C(X\times\T)$.
\epr

\bprf
Since $X$ is compact and $\Omega$ is closed, $X_1=X+\Omega$ becomes closed. 
The same reason shows that $X_2$ is closed.
It is easy to see that both $X_1$ and $X_2$ satisfy the conditions of $\omega$-invariance.
Note that $X_1\setminus X_2$ is a disjoint union of compact sets $X,X+\omega_1,\ldots,X+(K-1)\omega_1$.
Since $fS_i=S_i\sigma_{\omega_i}f=0$ for $i\neq 1$ and $f\in C(X_1\setminus X_2)\subset I_{X_2}/I_{X_1}$, we have $S_1 fS_1^*=\sigma_{-\omega_1}(f)S_1 S_1^*=\sigma_{-\omega_1}f-\sigma_{-\omega_1}f\sum_{i=2}^nS_i S_i^*=\sigma_{-\omega_1}f$. Therefore $I_{X_2}/I_{X_1}=\cspa\{S_\mu fS_\nu^*\mid \mu,\nu\in\W_n, f\in C(X)\}$.
For $(k,\mu),(l,\nu)\in \Z/K\Z\times\W_n^+$, let us define $e_{(k,\mu),(l,\nu)}=S_\mu S_1^k\chi {S_1^*}^lS_\nu^*\in I_{X_2}/I_{X_1}$ where $\chi$ is a characteristic function of $X$.
Then $\{e_{(k,\mu),(l,\nu)}\}$ satisfies the relation of matrix units and $\sum_{(k,\mu)\in \Z/K\Z\times\W_n^+}e_{(k,\mu),(k,\mu)}=1$ (strictly).
Since $e_{(0,\emptyset),(0,\emptyset)}=\chi$, we have $I_{X_2}/I_{X_1}\cong\K\otimes B$ where $B=\chi(I_{X_2}/I_{X_1})\chi$.
We have 
$$B=\cspa\{\chi S_\mu fS_\nu^* \chi\mid \mu,\nu\in\W_n, f\in C(X)\}=\cspa\{(S_1^K)^m f\mid m\in\Z, f\in C(X)\}.$$
Since $B$ is generated by $C(X)$ and a unitary $S_1^K\chi$ which commute with each other and since $B$ is globally invariant under the gauge action, we have $B\cong C(X\times\T)$.
Therefore we get $I_{X_2}/I_{X_1}\cong\K\otimes C(X\times\T)$.
\eprf

Let us choose $\gamma_0\in\Gamma$ and fix it.
Set $X_1=\gamma_0+\Omega$ and $X_2=\gamma_0+(\Omega\setminus\{0,\omega_1,\ldots,(K-1)\omega_1\}$ which is an $\omega$-invariant subset of $\Gamma$ by Corollary \ref{disc2}.
Since $[X_1]\times\T\supset Y_{([\gamma_0],\theta_0)}\supset [X_2]\times\T$, we have $I_{X_1}\subset P_{([\gamma_0],\theta_0)}\subset I_{X_2}$ for any $\theta_0\in\T$.
Taking $X=\{\gamma_0\}$ in Proposition \ref{local}, we get an isomorphism $I_{X_2}/I_{X_1}\cong\K\otimes C(\T)$ which sends $S_1^K\chi$ to $p\otimes z$ where $\chi\in C_0(X_1\setminus X_2)$ is a characteristic function of $\gamma_0$, $p\in\K$ is a minimal projection corresponding to $\chi$, and $z$ is the canonical generator of $C(\T)$.

\blem\label{P}
Under the isomorphism $I_{X_2}/I_{X_1}\cong\K\otimes C(\T)$ above, 
we have $P_{([\gamma_0],\theta_0)}/I_{X_1}\cong\K\otimes C_0(\T\setminus\{\theta_0\})$ for any $\theta_0\in\T$.
\elem

\bprf
Since $P_{([\gamma_0],\theta_0)}/I_{X_1}$ is an ideal of $I_{X_2}/I_{X_1}$,
all we have to do is to show $(P_{([\gamma_0],\theta_0)}/I_{X_1})\cap C^*(S_1^K\chi)\cong p\otimes C_0(\T\setminus\{\theta_0\})$.
For $\theta\in\T$, the map $\psi_{\gamma_0,\theta}:A\to\M_K$ vanishes on $A\cap I_{X_1}$.
Hence we can define the map $\psi'_{\theta}:A/(A\cap I_{X_1})\to\M_K$ so that $\psi'_{\theta}\circ\pi'=\psi_{\gamma_0,\theta}$ where $\pi'$ is the canonical surjection from $A$ to $A/(A\cap I_{X_1})$.
The image of $C^*(S_1^K\chi)\subset A/(A\cap I_{X_1})$ under $\psi'_{\theta}$ is contained in $\C 1\subset\M_K$.
One can show that the map $\psi'_{\theta}:C^*(S_1^K\chi)\to\C1$ is isomorphic to the evaluation map at $\theta\in\T$ from $p\otimes C_0(\T)\subset \K\otimes C_0(\T)$ under the isomorphism $C^*(S_1^K\chi)\cong p\otimes C_0(\T)$.
Noting that $(P_{([\gamma_0],\theta_0)}/I_{X_1})\cap C^*(S_1^K\chi)=(J_{Y_{([\gamma_0],\theta_0)}}/I_{X_1})\cap C^*(S_1^K\chi)$ by Lemma \ref{Y_I_Y1}, 
we get the desired isomorphism $(P_{([\gamma_0],\theta_0)}/I_{X_1})\cap C^*(S_1^K\chi)\cong C_0(\T\setminus\{\theta_0\})$.
\eprf

To prove $P_{([\gamma_0],\theta_0)}$ is primitive, we need the following observation which is inspired by \cite{aHR}.
Let $H=l^2((\Z/K\Z)\times\W_n^+)$ be a Hilbert space whose complete orthonormal system is given by $\{\xi_{k,\mu}\mid k\in \Z/K\Z, \mu\in\W_n^+\}$.
For $i=1,2,\ldots,n$, let us define $T_i\in B(H)$ by 
$$T_i(\xi_{k,\mu})=\left\{\begin{array}{ll}\xi_{k+1,\emptyset}&(\mbox{if }i=1,\mu=\emptyset),\\
\xi_{k,i\mu}&(\mbox{otherwise}).
\end{array}\right.$$
For $\gamma\in\Omega$, let us denote by $Q_\gamma\in B(H)$ a projection onto $\cspa\{\xi_{k,\mu}\mid k\omega_1+\omega_{\mu}=\gamma\}$.
One can easily see the following.

\blem\label{T&Q}
\benu
\item For $i=1,2,\ldots,n$, $T_i^*T_i=1$ and $\sum_{i=1}^nT_iT_i^*=1$.
\item $\sum_{\gamma\in\Omega}Q_\gamma=1$ $($strongly$)$.
\item For $i=1,2,\ldots,n$ and $\gamma\in\Omega$, $Q_\gamma T_i=\left\{\begin{array}{ll}T_i Q_{\gamma-\omega_i}&(\mbox{if }\gamma-\omega_i\in\Omega),\\
0&(\mbox{otherwise}).
\end{array}\right.$
\eenu
\elem

By this lemma, there exists a unique $*$-homomorphism $\varphi_1:\cp\to B(H)$ with $\varphi_1(S_i)=T_i$ and $\varphi_1(f)=\sum_{\gamma\in\Omega}f(\gamma_0+\gamma)Q_\gamma$ for $f\in C_0(\Gamma)$.

\blem\label{irrep}
We have $\varphi_1(I_{X_1})=0$ and $\varphi_1(I_{X_2})=K(H)$.
\elem

\bprf
Since $\varphi_1(f)=0$ for any $f\in C_0(\Gamma\setminus X_1)$, we have $\varphi_1(I_{X_1})=0$.
From $I_{X_2}=\cspa\{S_\mu f S_\nu^*\mid \mu,\nu\in\W_n, f\in C_0(\Gamma\setminus X_2)\}$ and $Q_{k\omega_1}=T_1^kQ_0{T_1^*}^k$, we get 
$$\varphi_1(I_{X_2})=\cspa\{T_\mu Q_{k\omega_1} T_\nu^*\mid \mu,\nu\in\W_n, k\in\Z/K\Z\}=\cspa\{T_\mu Q_0 T_\nu^*\mid \mu,\nu\in\W_n\}.$$
Writing $T_\mu=T_{\mu'}T_1^{l}$ and $T_\nu=T_{\nu'}T_1^{m}$ where $\mu',\nu'\in\W_n^+$ and $l,m\in\N$, we see that $T_\mu Q_0 T_\nu^*$ is a one rank operator from $\xi_{m',\nu'}$ to $\xi_{l',\mu'}$ where $m',l'\in\Z/K\Z$ are images of $m,l\in\N$ respectively.
Therefore $\varphi_1(I_{X_2})=K(H)$.
\eprf

Since $\varphi_1(\cp)\supset K(H)$, $\varphi_1$ is an irreducible representation.
Hence $\ker\varphi_1$ is a primitive ideal.
We will prove that $\ker\varphi_1=P_{([\gamma_0],\theta_0)}$ for some $\theta_0\in\T$, 

For $t\in\T$, let us define a unitary $U_t\in B(H)$ by $U_t(\xi_{k,\mu})=t^{|\mu|}\xi_{k,\mu}$.
One can easily see the following.

\blem\label{u}
\benu
\item $U_tQ_\gamma U_t^*=Q_\gamma$ for $\gamma\in\Omega$.
\item $U_tT_iU_t^*=tT_i$ for $i=2,3,\ldots,n$.
\item $U_tT_1U_t^*=tT_1+(1-t)V$ where $V\in B(H)$ is defined by $V(\xi_{k,\mu})=\delta_{\mu,\emptyset}\xi_{k+1,\emptyset}$.
\eenu
\elem

For $t\in\T$, let us define a $*$-automorphism $\beta'_t$ of $B(H)$ by $\beta'_t(x)=U_txU_t^*$ for $x\in B(H)$.
Since $\beta'_t(K(H))=K(H)$ for any $t\in\T$, we can extend the $*$-automorphism $\beta'_t$ of $B(H)$ to one of $B(H)/K(H)$ which is also denoted by $\beta'_t$.
Let us denote by $\varphi_2:\cp\to B(H)/K(H)$ the composition of $\varphi_1$ and a natural surjection $\pi:B(H)\to B(H)/K(H)$.

\blem\label{commute}
For any $t\in\T$, we have $\beta'_t\circ\varphi_2=\varphi_2\circ\beta_t$ where $\beta$ is the gauge action on $\cp$.
\elem

\bprf
The only non-trivial part is $\beta'_t(\pi(T_1))=t\pi(T_1)$ for $t\in\T$ which follows from the fact that $V$ in Lemma \ref{u} (iii) is a compact operator.
\eprf

\blem\label{kervarphi_2}
We have $\ker\varphi_2=I_{X_2}$.
\elem

\bprf
By Lemma \ref{commute}, $\ker\varphi_2$ is a gauge invariant ideal.
For $\gamma\in\Omega$, $Q_{\gamma}$ becomes a compact operator if and only if $\gamma=k\omega_1$ for some $k\in\Z/K\Z$.
Hence $X_{\ker\varphi_2}=\gamma_0+(\Omega\setminus\{0,\omega_1,\ldots,(K-1)\omega_1\})=X_2$.
Therefore we have $\ker\varphi_2=I_{X_2}$.
\eprf

\bpr\label{Pispri}
For any $\theta\in\T$, $P_{([\gamma_0],\theta)}$ is a primitive ideal.
\epr

\bprf
By Lemma \ref{irrep} and Lemma \ref{kervarphi_2}, we have $I_{X_1}\subset \ker\varphi_1\subset I_{X_2}$.
Since $\ker\varphi_1$ is primitive, the ideal $\ker\varphi_1/I_{X_1}$ of $I_{X_2}/I_{X_1}$ is also primitive.
Hence we have $\ker\varphi_1/I_{X_1}\cong\K\otimes C_0(\T\setminus\{\theta_0\})$ for some $\theta_0\in\T$.
By Lemma \ref{P}, we have $\ker\varphi_1=P_{([\gamma_0],\theta_0)}$.
Thus $P_{([\gamma_0],\theta_0)}$ is a primitive ideal.
For an arbitrary $\theta\in\T$, there exists $t\in\T$ with $\theta=t^K\theta_0$.
Hence $P_{([\gamma_0],\theta)}=\beta_t(P_{([\gamma_0],\theta_0)})$ is also primitive.
\eprf

In fact, we can prove $\ker\varphi_1=P_{([\gamma_0],1)}$, but we do not need it.

\bpr\label{Primsurj}
For $\gamma_0\in\Gamma$, the set of all primitive ideals $P$ satisfying $X_P=\gamma_0+\Omega$ is $\{P_{([\gamma_0],\theta)}\mid\theta\in\T\}$.
\epr

\bprf
By Proposition \ref{Pispri}, the ideal $P=P_{([\gamma_0],\theta)}$ is primitive with $X_P=\gamma_0+\Omega$ for any $\theta\in\T$.
Let $P$ be a primitive ideal of $\cp$ with $X_P=\gamma_0+\Omega$ for some $\gamma_0\in\Gamma$.
Then, $I_{X_1}\subset P$ and $I_{X_2}\not\subset P$ where $X_1=\gamma_0+\Omega$ and $X_2=\gamma_0+(\Omega\setminus\{0,\omega_1,\ldots,(K-1)\omega_1\})$.
The set of all primitive ideals $P$ satisfying $I_{X_1}\subset P$ and $I_{X_2}\not\subset P$ corresponds to the set of primitive ideals of $I_{X_2}/I_{X_1}\cong\K\otimes C(\T)$ bijectively (see, for example, \cite{Dix}).
Hence there is no primitive ideal $P$ satisfying $X_P=\gamma_0+\Omega$ other than $\{P_{([\gamma_0],\theta)}\mid\theta\in\T\}$.
\eprf

Now we can describe the primitive ideal space of $\cp$.

\blem\label{Priminj}
Let $\gamma_1,\gamma_2\in\Gamma$ and $\I_1,\I_2$ be non-empty sets of $\{1,2,\ldots,n\}$.
Then $I_{\gamma_1+\Omega_{\I_1}}=I_{\gamma_2+\Omega_{\I_2}}$ if and only if $\Omega_{\I_1}=\Omega_{\I_2}$ and $\gamma_1-\gamma_2\in\Omega_{\I_1}\cap (-\Omega_{\I_1})$.
\elem

\bprf
Obviously $I_{\gamma_1+\Omega_{\I_1}}=I_{\gamma_2+\Omega_{\I_2}}$ is equivalent to $\gamma_1+\Omega_{\I_1}=\gamma_2+\Omega_{\I_2}$.
If $\Omega_{\I_1}=\Omega_{\I_2}$ and $\gamma_1-\gamma_2\in\Omega_{\I_1}\cap (-\Omega_{\I_1})$, then $\gamma_1+\Omega_{\I_1}=\gamma_2+\Omega_{\I_2}$.
Conversely assume $\gamma_1+\Omega_{\I_1}=\gamma_2+\Omega_{\I_2}$ and denote it by $X$.
Then we have $\Omega_{\I_1}=\Omega_{\I_2}$ because $\Omega_{\I_j}=\{\gamma\in\Gamma\mid X+\gamma\subset X\}$ for $j=1,2$.
Hence we get $\gamma_1-\gamma_2\in\Omega_{\I_1}\cap (-\Omega_{\I_1})$.
The proof is completed.
\eprf

For non-empty sets $\I_1,\I_2$ of $\{1,2,\ldots,n\}$, we write $\I_1\sim\I_2$ if $\Omega_{\I_1}=\Omega_{\I_2}$.
Let us choose and fix representative elements of each equivalence classes of $\sim$ and denote by ${\mathcal I}$ the set of them.
Note that $\{1\}\in{\mathcal I}$ since $\I\sim\{1\}$ if and only if $\I=\{1\}$.
For each $\I\in{\mathcal I}$, we define a topological space $\Gamma_{\I}$ by $\Gamma_{\I}=\Gamma/(\Omega_{\I}\cap (-\Omega_{\I}))$ if $\I\neq\{1\}$ and $\Gamma_{\{1\}}=\Gamma'\times\T$.
For $[\gamma]\in\Gamma_{\I}$ with $\I\neq\{1\}$, we define a primitive ideal $P_{[\gamma]}$ by $I_{\gamma+\Omega_{\I}}$.
Note that if $[\gamma]=[\gamma']$ in $\Gamma_{\I}$, then $I_{\gamma+\Omega_{\I}}=I_{\gamma'+\Omega_{\I}}$.
For $([\gamma],\theta)\in \Gamma_{\{1\}}=\Gamma'\times\T$, the ideal $P_{([\gamma],\theta)}$ is defined in Definition \ref{DefP}.

\bpr\label{primitive2}
The map $\coprod_{\I\in{\mathcal I}}\Gamma_{\I}\ni y \mapsto P_y \in\Prim(\cp)$ is bijective where $\Prim(\cp)$ is the primitive ideal space of $\cp$.
\epr

\bprf
The map is injective by Lemma \ref{Priminj} and surjective by Proposition \ref{Primsurj}.
\eprf

The primitive ideal space $\Prim(\cp)$ is a topological space whose closed sets are given by $\{P\in\Prim(\cp)\mid I\subset P\}$ for ideals $I$.
We will investigate which subset of $\coprod_{\I\in{\mathcal I}}\Gamma_{\I}$ corresponds to a closed subset of $\Prim(\cp)$.

\blem\label{Primlem1}
For any $\gamma_0\in\Gamma$, there exists a compact neighborhood $X$ of $\gamma_0$ such that $X\cap (X+\gamma)=\emptyset$ for any $\gamma\in\Omega\setminus\{0\}$.
\elem

\bprf
By Proposition \ref{disc}, there exists a neighborhood $U$ of $0$ with $\gamma\notin U$ for any $\gamma\in\Omega\setminus\{0\}$.
Take a compact neighborhood $V$ of $0$ with $V-V\subset U$.
Then $X=\gamma_0+V$ is a desired compact neighborhood of $\gamma_0$.
\eprf

\blem\label{Primlem2}
We have $X+\Omega\not\supset\gamma+\Omega_{\I}$ for any $\gamma\in\Gamma$, $\I\neq\{1\}$, and any compact set $X$ of $\Gamma$.
\elem

\bprf
To derive a contradiction, assume $X+\Omega\supset\gamma+\Omega_{\I}$ for some $\gamma\in\Gamma$, $\I\neq\{1\}$, and some compact set $X$ of $\Gamma$.
Take $i\in\I$ with $i\neq 1$.
For any $k\in\N$, the element $\gamma-k\omega_i$ is in $\gamma+\Omega_{\I}$.
Hence there exists $\gamma_k\in X$ and $\mu_k\in\W_n$ with $\gamma-k\omega_i=\gamma_k+\omega_{\mu_k}$ for any $k\in\N$.
Since $X$ is compact, there exists a subsequence $\{k_l\}_{l\in\N}$ of $\N$ 
so that $\gamma_{k_l}$ converges to some point.
Thus $\{\omega_{\mu_{k_l}}+k_l\omega_i\}_{l\in\N}$ becomes 
a convergent sequence.
By the same argument as in the proof of Proposition \ref{disc}, 
we can show that $-\omega_i\in\Omega=\Omega_{\{1\}}$.
This contradicts the assumption for $\omega$.
\eprf

\blem\label{Primlem3}
Let $Y$ be a subset of $\Gamma'\times\T$.
If for any $[\gamma]\in\Gamma'$, there exists a compact neighborhood $[X_{\gamma}]$ of $[\gamma]$ such that $Y\cap ([X_{\gamma}]\times\T)$ is closed set, then $Y$ is closed.
\elem

\bprf
Take a net $\{([\gamma_\lambda],\theta_\lambda)\}$ in $Y$ converges to $([\gamma],\theta)\in\Gamma'\times\T$.
Eventually $[\gamma_\lambda]\in [X_{\gamma}]$  because $[X_\gamma]$ is a neighborhood of $[\gamma]$.
Then $([\gamma],\theta)\in Y$ since $Y\cap ([X_{\gamma}]\times\T)$ is closed.
Thus $Y$ is closed.
\eprf

\blem\label{Primlem4}
For any $\omega$-invariant subset $X$ of $\Gamma$, we have $I_{X}=\bigcap_{y\in [X]\times\T}P_y$.
\elem

\bprf
By Proposition \ref{rotation}, the ideal $I=\bigcap_{y\in [X]\times\T}P_y$ is gauge invariant.
Hence $I=I_X$ because $C_0(\Gamma)\cap I=\bigcap_{y\in [X]\times\T}(C_0(\Gamma)\cap P_y)=\bigcap_{\gamma\in X}C_0(\Gamma\setminus (\gamma+\Omega))=C_0(\Gamma\setminus X)$.
\eprf

\bpr\label{closed}
Let $Y$ be a subset of $\coprod_{\I\in{\mathcal I}}\Gamma_{\I}$ and set $Y_\I=Y\cap\Gamma_{\I}$ for $\I\in{\mathcal I}$.
The set $P_Y=\{P_y\mid y\in Y\}$ is closed in $\Prim(\cp)$ if and only if $Y_{\{1\}}$ is an $\omega$-invariant set of $\Gamma_{\{1\}}=\Gamma'\times\T$ and $Y_{\I}=\{[\gamma]\in\Gamma_{\I}\mid [\gamma+\Omega_{\I}]\times\T\subset Y_{\{1\}}\}$ for any $\I\in{\mathcal I}$ with $\I\neq\{1\}$.
\epr

\bprf
Let us take a subset $Y=\coprod_{\I\in{\mathcal I}}Y_\I$ of $\coprod_{\I\in{\mathcal I}}\Gamma_{\I}$.
If $Y_{\{1\}}$ is an $\omega$-invariant subset of $\Gamma_{\{1\}}=\Gamma'\times\T$, then we can define the ideal $I_{Y_{\{1\}}}$.
One can easily see that $\{([\gamma],\theta)\in\Gamma_{\{1\}}\mid I_{Y_{\{1\}}}\subset P_{([\gamma],\theta)}\}=Y_{\{1\}}$ and that for $\I\neq\{1\}$, $I_{Y_{\{1\}}}\subset I_{\gamma+\Omega_{\I}}$ if and only if $[\gamma+\Omega_{\I}]\times\T\subset Y_{\{1\}}$.
Therefore if $Y$ satisfies the condition in the statement, then $P_Y$ is closed in $\Prim(\cp)$.

Conversely, assume $P_Y$ is closed, i.e. there exists an ideal $I$ of $\cp$ so that $Y=\{y\in\coprod_{\I\in{\mathcal I}}\Gamma_{\I}\mid I\subset P_y\}$.
We first show that $Y_{\{1\}}$ is $\omega$-invariant.
Take $\gamma_0\in\Gamma$ arbitrarily.
By Lemma \ref{Primlem1}, there exists a compact neighborhood $X$ of $\gamma_0$ such that $X\cap (X+\gamma)=\emptyset$ for any $\gamma\in\Omega\setminus\{0\}$.
If we set $X_1=X+\Omega$ and $X_2=X+(\Omega\setminus\{0,\omega_1,\ldots,(K-1)\omega_1\})$, then $I_{X_2}/I_{X_1}\cong\K\otimes C(X\times\T)$ by Proposition \ref{local}.
The subset $\{P\in\Prim(\cp)\mid I_{X_1}\subset P, I_{X_2}\not\subset P\}$ of $\Prim(\cp)$ is homeomorphic to $\Prim(I_{X_2}/I_{X_1})\cong X\times\T$.
By Lemma \ref{Primlem2}, $X_1\not\supset\gamma+\Omega_{\I}$ for any $\gamma\in\Gamma$ and for any $\I\neq\{1\}$.
Hence 
$$\big\{x\in \coprod_{\I\in{\mathcal I}}\Gamma_{\I}\mid I_{X_1}\subset P_x, I_{X_2}\not\subset P_x\big\}=[X]\times\T\subset \Gamma_{\{1\}}.$$
Therefore $[X]\times\T\ni x\mapsto P_x\in\Prim(\cp)$ is a homeomorphism from $[X]\times\T$ whose topology is the relative topology of $\Gamma'\times\T$ to the subset $\{P\in\Prim(\cp)\mid I_{X_1}\subset P, I_{X_2}\not\subset P\}$ of $\Prim(\cp)$ (note that $X$ is homeomorphic to $[X]$).
The set $Y\cap ([X]\times\T)\subset \Gamma_{\{1\}}$ is closed in $[X]\times\T$ because $P_Y$ is closed.
By Lemma \ref{Primlem3}, the subset $Y_{\{1\}}$ is closed in $\Gamma_{\{1\}}$.
If $([\gamma],\theta)\in Y_{\{1\}}$, then $([\gamma+\omega_i],\theta')\in Y_{\{1\}}$ for any $i\in\{2,3,\ldots,n\}$ and $\theta'\in\T$ because $P_{([\gamma],\theta)}\subset P_{([\gamma+\omega_i],\theta')}$.
Therefore $Y_{\{1\}}$ is an $\omega$-invariant subset of $\Gamma_{\{1\}}$.
Take $\I\in{\mathcal I}$ with $\I\neq\{1\}$ and $[\gamma]\in\Gamma_\I$.
Since $I_{\gamma+\Omega_\I}=\bigcap_{y\in [\gamma+\Omega_\I]\times\T}P_y$ by Lemma \ref{Primlem4}, the element $[\gamma]$ is in $Y_\I$ if and only if $[\gamma+\Omega_\I]\times\T\subset Y_{\{1\}}$.
Therefore $Y$ satisfies the condition in the statement.
\eprf

By the proposition above, we get the following.

\bthm\label{idestr2}
When $\omega$ does not satisfy Condition \ref{cond},
there is a one-to-one correspondence between the set of ideals of $\cp$ and the set of $\omega$-invariant subsets of $\Gamma'\times\T$.
Hence for any ideal $I$ of $\cp$, we have $I=I_{Y_I}$. 
\ethm

\bprf
There is a one-to-one correspondence between the set of ideals of $\cp$ and the closed subset of $\Prim(\cp)$.
By Proposition \ref{closed}, the closed subset of $\Prim(\cp)$ corresponds bijectively to the set of $\omega$-invariant subsets of $\Gamma'\times\T$.
\eprf

\section{The strong Connes spectrum and the K-groups of $\cp$}

As a consequence of knowing all ideals of $\cp$, we can compute the strong Connes spectrum of the action $\alpha^{\omega}:G\curvearrowright\On$.
We recall the definition of the strong Connes spectrum.

\bde
Let $\alpha:G\curvearrowright A$ be an action of an abelian group $G$, whose dual group is $\Gamma$, on a $C^*$-algebra $A$.
The strong Connes spectrum $\widetilde{\Gamma}(\alpha)$ of $\alpha$ is defined by 
$$\widetilde{\Gamma}(\alpha)=\{\gamma\in\Gamma\mid \widehat{\alpha}_\gamma(I)\subset I, \mbox{for any ideal } I \mbox{ of } A\rtimes_\alpha G\},$$
where $\widehat{\alpha}:\Gamma\curvearrowright A\rtimes_\alpha G$ is the dual action of $\alpha$.
\ede

For each action $\alpha$, the strong Connes spectrum $\widetilde{\Gamma}(\alpha)$ is a closed subsemigroup of $\Gamma$.
We remark that in the original paper \cite{Ki}, A. Kishimoto defined the strong Connes spectrum in a different way and proved that his definition is equivalent to the definition above (see, \cite[Lemma 3.4]{Ki}).

In our setting, the dual actions $\widehat{\alpha^\omega}:\Gamma\curvearrowright \cp$ are characterized by $\widehat{\alpha^\omega}_\gamma(S_\mu fS_\nu^*)=S_\mu \sigma_{\gamma}fS_\nu^*$ for $\mu,\nu\in\W_n$, $f\in C_0(\Gamma)$ and $\gamma\in\Gamma$.

\bpr\label{SCS}
Let $\omega$ be an element of $\Gamma^n$.
The strong Connes spectrum $\widetilde{\Gamma}(\alpha^{\omega})$ of the action $\alpha^{\omega}$ is $\bigcap_{i=1}^n\Omega_{\{i\}}$.
\epr

\bprf
First we consider the case that $\omega$ satisfies Condition \ref{cond}.
Since the correspondence between ideals of $\cp$ and $\omega$-invariant subsets of $\Gamma$ is one-to-one by Theorem \ref{idestr1}, 
$\widehat{\alpha^\omega}_\gamma(I)\subset I$ if and only if $X_I-\gamma\supset X_I$ for an ideal $I$ and $\gamma\in\Gamma$.
For any $i\in\{1,2,\ldots,n\}$, the set $\Omega_{\{i\}}$ is an $\omega$-invariant set satisfying $\{\gamma\in\Gamma\mid \Omega_{\{i\}}+\gamma\subset\Omega_{\{i\}}\}=\Omega_{\{i\}}$.
Therefore $\widetilde{\Gamma}(\alpha^{\omega})\subset\bigcap_{i=1}^n\Omega_{\{i\}}$.
We have $X\supset X+\bigcap_{i=1}^n\Omega_{\{i\}}$ for any $\omega$-invariant set $X$ because for any $x\in X$ there exists $i$ with $x+\Omega_{\{i\}}\subset X$.
we have $\widetilde{\Gamma}(\alpha^{\omega})\supset\bigcap_{i=1}^n\Omega_{\{i\}}$.
Thus $\widetilde{\Gamma}(\alpha^{\omega})=\bigcap_{i=1}^n\Omega_{\{i\}}$ in the case that $\omega$ satisfies Condition \ref{cond}.

Next we consider the case that $\omega$ does not satisfy Condition \ref{cond}.
In this case, the set $\bigcap_{i=1}^n\Omega_{\{i\}}$ coincides with $\Omega=\{\omega_\mu\mid\mu\in\W_n\}$.
Since $\Omega$ is an $\omega$-invariant subset of $\Gamma$ and $\{\gamma\in\Gamma\mid \widehat{\alpha^\omega}_\gamma(I_\Omega)\subset I_\Omega\}=\{\gamma\in\Gamma\mid \Omega+\gamma\subset\Omega\}=\Omega$, we have $\widetilde{\Gamma}(\alpha^{\omega})\subset\Omega$.
For any $\omega$-invariant subset $Y$ of $\Gamma'\times\T$, we have $([\gamma+\omega_\mu],\theta)\in Y$ for any $\mu\in\W_n$ and any $([\gamma],\theta)\in Y$.
Since the correspondence between ideals of $\cp$ and $\omega$-invariant subsets of $\Gamma'\times\T$ is one-to-one by Theorem \ref{idestr2}, 
we have $\widehat{\alpha^\omega}_\gamma(I)\subset I$ for any ideal $I$ of $\cp$ and for any $\gamma\in\Omega$.
Hence $\widetilde{\Gamma}(\alpha^{\omega})\supset\Omega$.
Therefore also in the case that $\omega$ does not satisfy Condition \ref{cond}, we have $\widetilde{\Gamma}(\alpha^{\omega})=\bigcap_{i=1}^n\Omega_{\{i\}}$.
\eprf

\bre
The inclusion $\widetilde{\Gamma}(\alpha^{\omega})\subset\bigcap_{i=1}^n\Omega_{\{i\}}$ had been already proved by A. Kishimoto \cite{Ki}.
\ere

The crossed product $\cp$ is a Cuntz-Pimsner algebra.
Let $E=C_0(\Gamma)^n$ be a right $C_0(\Gamma)$ module.
The left $C_0(\Gamma)$ module structure of $E$ is given by 
$$f\cdot(f_1,f_2,\ldots,f_n)=(\sigma_{\omega_1}(f)f_1,\sigma_{\omega_2}(f)f_2,\ldots,\sigma_{\omega_n}(f)f_n)\in E$$
for $f\in C_0(\Gamma)$ and $(f_1,f_2,\ldots,f_n)\in E$.

\bpr\label{CP}
The crossed product $\cp$ is isomorphic to the Cuntz-Pimsner algebra ${\mathcal O}_E$.
\epr

\bprf
The inclusion $C_0(\Gamma)\hookrightarrow\cp$ and $E\ni (f_1,f_2,\ldots,f_n)\mapsto \sum_{i=1}^n S_if_i\in\cp$ satisfies the conditions in \cite[Theorem 3.12]{Pi} (for example, the condition (4) is equivalent to saying that $\sum_{i=1}^nS_i\sigma_{\omega_i}(f)S_i=f$ for any $f\in C_0(\Gamma)$). 
Hence there exists a $*$-homomorphism $\varphi:{\mathcal O}_E\to\cp$ which is surjective since $\cp$ is generated by $\{\sum_{i=1}^n S_if_i\mid f_i\in C_0(\Gamma)\}$.
One can show that $\varphi$ is injective by using Proposition \ref{isom}.
Thus $\cp$ is isomorphic to ${\mathcal O}_E$.
\eprf

The ideal structures of Cuntz-Pimsner algebras were investigated in \cite{KPW} when Hilbert bimodules are finitely generated.
Our Hilbert bimodule $E$ is finitely generated if and only if the group $G$ is discrete.
When $G$ is discrete, we can know the detailed structure of $\cp$ without using the result in \cite{KPW} (see subsection \ref{G=disc}).
Thanks to considering our algebra $\cp$ as a Cuntz-Pimsner algebra,
we can compute the K-groups of it by \cite[Theorem 4.8]{Pi}.

\bpr
Let $\omega$ be an element of $\Gamma^n$.
The following sequence is exact:
$$\begin{CD}
K_0(C_0(\Gamma)) @>{\rm id}-\sum_{i=1}^n(\sigma_{\omega_i})_*>> K_0(C_0(\Gamma)) @>\iota_*>>  K_0(\cp) \\
@AAA @. @VVV \\
K_1(\cp) @<\iota_*<< K_1(C_0(\Gamma)) @<{\rm id}-\sum_{i=1}^n(\sigma_{\omega_i})_*<< K_1(C_0(\Gamma)),
\end{CD}$$
where $\iota$ is the embedding $\iota:C_0(\Gamma)\hookrightarrow\cp$.
\epr

\bprf
Let us denote by ${\mathcal T}_n$ the Cuntz-Toeplitz algebra, which is generated by $n$ isometries $T_1,T_2,\ldots,T_n$ satisfying $\sum_{i=1}^nT_iT_i^*<1$.
There is a surjection $\pi:{\mathcal T}_n\to\On$ with $\pi(T_i)=S_i$ for $i=1,2,\ldots,n$.
The kernel of $\pi$ is isomorphic to $\K$.
If we define an action $\overline{\alpha}^\omega:G\curvearrowright{\mathcal T}_n$ by $\overline{\alpha}^\omega_t(T_i)=\ip{t}{\omega_i}T_i$ for $t\in G$ and $i=1,2,\ldots,n$, then the kernel of $\pi$ is invariant under this action and $\pi\circ\overline{\alpha}^\omega_t={\alpha}^\omega_t\circ\pi$ for any $t\in G$.
Hence there exists a short exact sequence $0\to\K\rtimes_{\overline{\alpha}^\omega}G\to{\mathcal T}_n\rtimes_{\overline{\alpha}^\omega}G\to\cp\to 0$.
One can see that ${\mathcal T}_n\rtimes_{\overline{\alpha}^\omega}G$ is isomorphic 
to ${\mathcal T}_E$ in a similar way to Proposition \ref{CP}.
The $C^*$-algebra $\K\rtimes_{\overline{\alpha}^\omega}G$ is isomorphic to $\K\otimes C_0(\Gamma)$.
The subalgebra $\C1\rtimes_{\overline{\alpha}^\omega}G$ of ${\mathcal T}_n\rtimes_{\overline{\alpha}^\omega}G$ is isomorphic to $C_0(\Gamma)$.
The inclusion $C_0(\Gamma)\hookrightarrow{\mathcal T}_n\rtimes_{\overline{\alpha}^\omega}G$ induces a KK-equivalence between $C_0(\Gamma)$ and ${\mathcal T}_n\rtimes_{\overline{\alpha}^\omega}G$ 
whose inverse is given by a Kasparov bimodule $({\mathcal E}_+\oplus{\mathcal E}_+,\pi_0\oplus\pi_1,T)\in KK({\mathcal T}_n\rtimes_{\overline{\alpha}^\omega}G,C_0(\Gamma))$ where ${\mathcal E}_+=\bigoplus_{k=0}^\infty E^{\otimes k}$ is a right $C_0(\Gamma)$ module, $\pi_0:{\mathcal T}_n\rtimes_{\overline{\alpha}^\omega}G\to L({\mathcal E}_+)$ is the natural representation, $\pi_1:{\mathcal T}_n\rtimes_{\overline{\alpha}^\omega}G\to L(\bigoplus_{k=1}^\infty E^{\otimes k})\subset L({\mathcal E}_+)$ is the representation obtained from the universal property of ${\mathcal T}_n\rtimes_{\overline{\alpha}^\omega}G$ and $T:{\mathcal E}_+\oplus{\mathcal E}_+\to{\mathcal E}_+\oplus{\mathcal E}_+$ is defined by $T(\xi\oplus\zeta)=\zeta\oplus\xi$ (for the detail, see Section 4 in \cite{Pi}).
To show that the 6-term exact sequence obtained from the short exact sequence $0\to\K\rtimes_{\overline{\alpha}^\omega}G\to{\mathcal T}_n\rtimes_{\overline{\alpha}^\omega}G\to\cp\to 0$ is the desired one, it suffices to see that the element $({\mathcal E}_+\oplus{\mathcal E}_+,(\varphi\circ\pi_0)\oplus(\varphi\circ\pi_1),T)\in KK(C_0(\Gamma),C_0(\Gamma))$ coincides with ${\rm id}-\sum_{i=1}^n(\sigma_{\omega_i})_*$ where $\varphi:C_0(\Gamma)\to {\mathcal T}_n\rtimes_{\overline{\alpha}^\omega}G$ is given by $\varphi(f)=(1-\sum_{i=1}^nT_iT_i^*)f$ (note that $1-\sum_{i=1}^nT_iT_i^*$ is a minimal projection of the kernel of $\pi$ which is isomorphic to $\K$).
A routine computation shows that $\varphi\circ\pi_0$ vanishes on $\bigoplus_{k=1}^\infty E^{\otimes k}$ and $\varphi\circ\pi_1$ vanishes on $\bigoplus_{k=2}^\infty E^{\otimes k}$ and on $E^{\otimes 0}$ ($=C_0(\Gamma)$).
Thus 
\begin{align*}
({\mathcal E}_+\oplus{\mathcal E}_+,(\varphi\circ\pi_0)\oplus(\varphi\circ\pi_1),T)&=(C_0(\Gamma)\oplus C_0(\Gamma),(\varphi\circ\pi_0)\oplus 0,T) + (E\oplus E, 0\oplus (\varphi\circ\pi_1),T)\\
&=(C_0(\Gamma)\oplus C_0(\Gamma),{\rm id}\oplus 0,T) + \sum_{i=1}^n (C_0(\Gamma)\oplus C_0(\Gamma), 0\oplus \sigma_{\omega_i},T)\\
&={\rm id}-\sum_{i=1}^n(\sigma_{\omega_i})_*.
\end{align*}
\eprf

\section{Examples and remarks}

\subsection{When $G$ is compact}

When $G$ is compact, its dual group $\Gamma$ becomes discrete.
In this case, for any $\omega\in\Gamma^n$ the crossed product $\cp$ becomes a graph algebra of some skew product graph which is row-finite (see \cite{KP}) and a part of our results here has been already proved in \cite{BPRS}.
There are many graph algebras which are not isomorphic to $\cp$, and it should be interesting to determine the ideal structures of such algebras.
Our technique here may help it.
We may consider our $C^*$-algebras $\cp$ as a continuous counterpart of graph algebras.
It seems to be interesting to define and examine graph algebras of continuous graph.

\subsection{When $G$ is discrete}\label{G=disc}

When $G$ is discrete, its dual group $\Gamma$ becomes compact.
Let us choose $\omega\in\Gamma^n$ and fix it.
Let us denote by $\overline{\Omega}$ a closed semigroup generated by $\omega_1,\omega_2,\ldots,\omega_n$.

\bpr
When $G$ is discrete, we have $-\omega_i\in\overline{\Omega}$ for $i=1,2,\ldots,n$.
\epr

\bprf
Let us take $i\in\{1,2,\ldots,n\}$.
Since $\Gamma$ is compact, a sequence $\{k\omega_i\}_{k=1}^\infty$ has a subsequence $\{k_l\omega_i\}_{l=1}^\infty$ which converges to some element in $\Gamma$.
For any $l$, we have $(k_{l+1}-k_l-1)\omega_i\in\Omega$ because $k_{l+1}>k_l$.
Hence $-\omega_i=\lim_{l\to\infty}(k_{l+1}-k_l-1)\omega_i\in\overline{\Omega}$.
\eprf

The following are easy consequences of above proposition.

\bco
Any $\omega\in\Gamma^n$ satisfies Condition \ref{cond}.
\eco

\bco
The set $\overline{\Omega}$ becomes a closed subgroup of $\Gamma$ and the set of all $\omega$-invariant subsets of $\Gamma$ is one-to-one correspondent to the set of all closed subset of $\Gamma/\overline{\Omega}$.
\eco

By two corollaries above, the set of all ideals of $\cp$ is one-to-one correspondence to the set of all closed subset of $\Gamma/\overline{\Omega}$.
In fact, we can examine the ideal structures of $\cp$ directly and more structures of the crossed product.
Let $G'$ be a quotient of $G$ by the closed subgroup
$$\{t\in G\mid \alpha_t^\omega=\mbox{id}\}=\{t\in G\mid \ip{t}{\omega_i}=0 \mbox{ for }i=1,2,\ldots,n\}=\{t\in G\mid \ip{t}{\gamma}=0 \mbox{ for any }\gamma\in\overline{\Omega}\}.$$
Then the dual group of $G'$ is naturally isomorphic to $\overline{\Omega}$.
Since $\omega\in\overline{\Omega}^n$, we can define an action $\alpha^\omega:G'\curvearrowright\On$.
The crossed product $\On{\rtimes_{\alpha^\omega}}G'$ becomes simple by Theorem \ref{simple}.
The crossed product $\cp$ becomes a continuous field over the space $\Gamma/\overline{\Omega}$ whose fiber of any point is isomorphic to the simple $C^*$-algebra $\On{\rtimes_{\alpha^\omega}}G'$ which is purely infinite (see \cite{KK2} or \cite{Ka}).

\subsection{When $G=\R$}

When $G$ is the real group $\R$, its dual group $\Gamma$ is also $\R$.
We define three types for elements of $\R^n$.

\bde
Let $\omega=(\omega_1,\omega_2,\ldots,\omega_n)\in\R^n$.
The element $\omega$ is said to be of type $(+)$ 
if all the $\omega_i$'s have the same sign, 
and to be of type $(-)$ if there exist $i,j$ such that $\omega_i<0<\omega_j$. 
Otherwise, the element $\omega$ is said to be of type $(0)$.
\ede

Namely $\omega$ is of type $(0)$ if and only if there exists 
$i\in\{1,2,\ldots,n\}$ such that $\omega_i=0$ and all the other 
$\omega_i$'s have the same sign.
When $\omega$ is of type $(+)$ or $(-)$, the set $\Omega_{\{i\}}$ coincides with the closed group generated by $\omega_1,\omega_2,\ldots,\omega_n$ for any $i=1,2,\ldots,n$.
An element $\omega\in\R^n$ is called aperiodic if the closed group generated by $\omega_1,\omega_2,\ldots,\omega_n$ is $\R$.
By Theorem \ref{simple}, we have the following.

\bpr
For $\omega\in\R^n$, the crossed product $\cpr$ is simple if and only if $\omega$ is aperiodic and of type $(+)$ or $(-)$.
\epr

When $\omega$ is of type $(+)$ or $(-)$ and not aperiodic, the crossed product $\cpr$ is isomorphic to a mapping torus whose fiber is the simple $C^*$-algebra $\On{\rtimes_{\alpha^{\omega'}}}\T$ where $\omega'=(\omega_1/K,\omega_2/K,\ldots,\omega_n/K)\in\Z^n$ and $K$ is the (positive) generator of the closed group generated by $\omega_1,\omega_2,\ldots,\omega_n$ which is isomorphic to $\Z$.
Hence in this case, $\Prim(\cpr)\cong\T$ and the set of ideals of $\cpr$ corresponds to the set of closed sets of $\T$.
The case that $\omega$ is of type $(0)$ is more complicated.
When $\omega$ is of type $(0)$, the set $\Omega=\{\omega_\mu\mid\mu\in\W_n\}$ is closed and a closed set $X\subset\R$ is $\omega$-invariant if and only if $X+\Omega\subset X$.
We can prove the proposition below in a similar way to the proof of 
Proposition \ref{local}.

\bpr
Let $\omega\in\R^n$ be of type $(0)$ and $X$ be an $\omega$-invariant set.
Set $X'=\bigcup_{\omega_i\neq 0} (X+\omega_i)$.
Any closed set $X_1$ with $X'\subset X_1\subset X$ is $\omega$-invariant, and $I_{X_1}/I_X\cong\K\otimes C(X\setminus X_1)\otimes{\mathcal O}_k$ where $k$ is the number of $i$ with $\omega_i=0$ and ${\mathcal O}_1=C(\T)$.
\epr

One can easily see that an element $\omega\in\R^n$ does not satisfy Condition \ref{cond} if and only if $\omega$ is of type $(0)$ and the number of $i$ with $\omega_i=0$ is $1$.

\bre
When $\omega$ is of type $(+)$, 
the crossed product $\cpr$ becomes stable and projectionless \cite{KK1}.
In the forthcoming paper \cite{Ka}, 
we will show that $\cpr$ is AF-embeddable in this case.
More generally, we will give one sufficient condition for 
crossed product $\cp$ becomes AF-embeddable in \cite{Ka}.
As a consequence of it, we will show that $\cp$ is either AF-embeddable 
or purely infinite when it is simple.
\ere

\end{document}